\RequirePackage{fix-cm}
\documentclass[smallextended]{svjour3}       
\usepackage[utf8]{inputenc}
\usepackage{graphicx}
\usepackage{color}
 \usepackage{mathptmx}      
\usepackage{latexsym}
\usepackage{amsmath,amssymb}
\begin{document}

\title{Nonnegative partial $s$-goodness for the equivalence of a  0-1 linear program to weighted linear programming
}


\author{Meijia Han         \and
        Wenxing Zhu 
}


\institute{Meijia Han \at
              Center for Discrete Mathematics and Theoretical Computer Science, Fuzhou University, Fuzhou 350116, China \\
              \email{1099091207@qq.com}           
           \and
           Wenxing Zhu \at
               Center for Discrete Mathematics and Theoretical Computer Science, Fuzhou University, Fuzhou 350116, China \\
              \email{wxzhu@fzu.edu.cn}           
}

\date{Received: date / Accepted: date}

\maketitle

\begin{abstract}
The 0-1 linear programming problem with nonnegative constraint matrix and objective vector $e$ origins from many NP-hard combinatorial optimization problems. In this paper, we consider recovering an optimal solution to the problem from a weighted linear programming. We first formulate the problem equivalently as a  sparse optimization problem. Next, we consider the consistency of the optimal solution of the sparse optimization problem and the  weighted linear programming problem. In order to achieve this, we establish nonnegative partial $s$-goodness of the constraint matrix and the weighted vector. Further, we use two quantities to characterize a sufficient condition and necessary condition for the nonnegative partial $s$-goodness. However, the two quantities are difficult to calculate, therefore, we provide a computable upper bound for one of the two quantities to verify the nonnegative partial $s$-goodness.
Finally, we give three examples to illustrate that our theory is effective and verifiable.

\end{abstract}

\section{Introduction}\label{introduction}
\label{intro}
In this paper, we  consider the integer programming (IP) with  linear inequalities:
\begin{equation}\label{IP}
\min_x \{\sum\limits_{i=1}^{n} x_i: ~Ax\geq b,~ x\in\{0, 1\}^n\},
\end{equation}
where $A\geq0\in R^{m\times n}$ is a given matrix, $b\in R^m$ is a given vector. Problem (\ref{IP}) arises from many combinatorial optimization problems, such as the minimum vertex covering \cite{Chen2016}, the maximum independent set, and the max-cut problems \cite{BILLIONNET2010}, etc. These problems have been proven to be NP-hard, and many  exact  and heuristics or approximation algorithms have been proposed, including cutting-plane method \cite{GOMORY1958}, branch and bound \cite{DAKIN1965} and local branch \cite{Lodi2010} algorithms, relaxation induced neighbourhood search \cite{Danna2005}, and objective scaling ensemble approach \cite{Zhang2020}.

It is worth to note that most of the above mentioned algorithms are based on the corresponding linear programming relaxation. Therefore, it is natural to ask under which conditions the integer programming problem \eqref{IP} is  equivalent to its corresponding linear programming relaxation problem:
\begin{eqnarray*}\nonumber
\min_x \{\sum\limits_{i=1}^{n} x_i: ~Ax\geq b,~ 0\leq x\leq1\}.
\end{eqnarray*}
At present, there are some theories for ensuring that the above linear programming problem has integral solution.

The first well known theory is  totally unimodular  ($TUM$ for short), under which it holds that

 \begin{theorem} (\cite{HOFFMAN1976})
If $A$ is totally unimodular and $b$ is an integer vector, then the vertices of the polyhedron $P = \{x:~Ax \leq b,~ 0\leq x\leq1\}$ are integral.
 \end{theorem}

Another well known theory is the totally dual integrality ($TDI$ for short) proposed by Edmonds and Giles \cite{Edmonds1977}, which is a weaker sufficient condition than TUM.
\begin{theorem}(\cite{Edmonds1977})
If the linear system $Ax \leq b$ is TDI and  b is integer valued, then $P = \{x: ~Ax \leq b,~ 0\leq x\leq1\}$ is an integral polyhedron.
\end{theorem}

Note that the linear programming  relaxation of problem \eqref{IP} can be written into  the form of linear complementarity problem (LCP):
\begin{equation}\label{lcp}
\begin{array}{l}
q+Mz\geq0,\\
z^T(q+Mz)=0,\\
z\geq0.
\end{array}
\end{equation}
where $z\in R^n$ is the vector of variables.
Assuming that the solution set of LCP is non-empty, Chandrasekaran et al. \cite{Chandrasekaran1998} gave  the class $I$ of integral matrices $M$ if the corresponding LCP has an integer solution for each integral vector $q$. 
Utilizing TDI, Dubey and Neogy \cite{Dubey2018}  obtained some new conditions for the existence of an integer solution to LCP with a hidden $Z-$matrix and hidden $K-$matrix. These results are extended from $TUM$ or $TDI$.

In this paper, we wish to establish a sufficient and necessary condition other than $TUM$ and $TDI$, such that the integer program \eqref{IP} is solvable by  linear programming relaxation.
We consider providing an adjustable weight $c$, $0<c\leq1$,
such that the optimal solution of the weighted linear programming relaxation problem
\begin{eqnarray}\label{wlp}
\min_x\{c^Tx:Ax\geq b, 0\leq x\leq1\},
\end{eqnarray}
is an optimal solution of problem  \eqref{IP}.

Let the $l_0-$norm $\|x\|_0$ be defined as the number of nonzero elements in the vector $x$. In Section \ref{chaper2} of this paper, we will prove that  problem \eqref{IP} can be written equivalently as the $l_0-$norm minimization problem
\begin{equation}\label{la1}
\min_x\{\|x\|_0:Ax\geq b,0\leq x\leq 1\}.
\end{equation}
Therefore, we expect to establish  a sufficient condition through the sparse optimization approach, such that problems \eqref{la1} and \eqref{wlp} have the same unique optimal solution, and thus provides an optimal solution to problem \eqref{IP}.

Eq. \eqref{la1} is a sparse optimization problem. In this field, the problem
\begin{equation}\label{wo}
\min\limits_x\{\|x\|_0: A_1x+A_2y=b\}
\end{equation}
has been studied mostly, for the optimal solution of $\min\limits_x\{\|x\|_1: A_1x+A_2y=b\}$ to be an optimal solution of the above problem. The proposed well known conditions are the partial restricted isometry property (PRIP) of the parameter $\delta^r_{s-r}$ \cite{Bandeira2013},  and the partial null space property (PNSP) of order $s-r$ \cite{Bandeira2013}. 
In \cite{Zhao2014,zhao2017,Zhao2018}, the authors provided the $k-$th order range space property ($k$-RSP) of  the constraint matrix $A$ and a given matrix $W$, under which the $l_0-$norm minimization problem has the same  unique optimal solution as the  problem  $\min\limits_x\{\|W_1x+W_2y\|_1: A_1x+A_2y=b\}$. The condition was further extended to the  $l_0-$norm minimization problem with non-negative constraints \cite{zhao2012}.

Another condition is the $s$-goodness of the constraint matrix $A$ \cite{Juditsky2011}, under which the unique optimal solution of the $l_1$-norm minimization problem  $\min\limits_x\{\|x+y\|_1: A_1x+A_2y=b\}$ is exactly an optimal one of problem \eqref{wo}. In \cite{Kong2014}, the condition was extended for considering  the problem
\begin{equation}
\min\limits_x\{\|x\|_1: A_1x+A_2y=b\},
\end{equation}
and the partial $s$-goodness of the constraint matrix $(A_1, A_2)$ was established, for the exact partial $s$-sparse optimal solution via the partial $l_1-$norm minimization.

In this paper, we will propose the nonnegative partial $s$-goodness condition for problems \eqref{wlp} and \eqref{la1}, such that they have the same unique optimal solution.  Specifically, we propose a definition of nonnegative partial $s$-goodness and its characterization $\gamma_{s,K}(\cdot)$ and $\hat{\gamma}_{s,K}(\cdot)$ with respect to the constraint matrices and the coefficients of the objective function. On this basis, we give  a computable upper bound of $\hat{\gamma}_{s,K}(\cdot)$, and thus obtain verifiable sufficient conditions for the nonnegative partial $s$-goodness. Through which we provide conditions such that the optimal solution of problem \eqref{IP} can be obtained from problem \eqref{wlp}.

This paper is organized as follows. Section \ref{chaper2} gives an example of problem \eqref{IP}, and proves the equivalence between problems (\ref{IP}) and (\ref{la1}). In Section \ref{chaper3}, we define nonnegative partial $s$-goodness and its characterization for a constraint matrix $A^\prime$ and the coefficient $c$ of the objective function in problem \eqref{wlp}. Moreover, we derive necessary condition and sufficient condition for $(A^\prime,c)$ to be nonnegative partial $s$-goodness,  and discuss efficiently computable upper bound of $\hat{\gamma}_{s,K}(\cdot)$ in Section \ref{chaper4}.
In Section \ref{chaper5}, we give a heuristic algorithm and three examples to show the feasibility of the proposed theory. Finally, Section \ref{chaper6} concludes this paper.

\section{Problem reformulation}\label{chaper2}
In this section, we prepare some preliminary work for the consequent research. First, we give an example of the considered problem \eqref{IP}. Then, we show that the considered integer programming problem \eqref{IP} can be converted equivalently to an $l_0$-norm minimization problem.

\subsection{An example of the considered integer programming problem}\label{chaper21}
We take the maximum independent set problem as an example to show that it is a special form of  problem \eqref{IP}.
Given an undirected graph $G=(V,E)$, where $V=\{1,\cdots,n\}$ is the set of vertices and $E$ is the set of edges,   the maximum independent set problem  can be formulated as
\begin{equation}\label{MIS}
\begin{array}{cl}
\max\limits_x & \sum\limits_{i=1}^{n}x_i\\
s.t. &Ax\leq 1\\
&x\in\{0,1\}^n,
\end{array}\end{equation}
where $A$ is the adjacency matrix of graph $G$.

Let $\tilde{x}_i=1-x_i$, $i=1, 2, \cdots, n$. Then  problem (\ref{MIS}) can be written equivalently as
\begin{equation*} 
\begin{array}{cl}
\min\limits_{\tilde{x}}&\sum\limits_{i=1}^{n}\tilde{x}_i\\
s.t. &A\tilde{x}\geq 1\\
&\tilde{x}\in\{0,1\}^n,
\end{array}\end{equation*}
which is a special form of   problem \eqref{IP}.

%

\subsection{Equivalence between integer programming problem and  $l_0$-norm minimization problem}\label{chaper22}
In this subsection, we  show the equivalence between  optimal solution of the $l_0$-norm minimization problem  \eqref{la1} and the integer programming problem \eqref{IP}  after a certain operation. First, we have the following result.



\begin{theorem}\label{thm21}
For any optimal solution $x^*$ of the $l_0$-norm minimization problem $\eqref{la1}$, $\hat{x}^*=(\lceil x_1^* \rceil, \lceil x_2^* \rceil,$ $\cdots, \lceil x_n^* \rceil)^T$ is an optimal solution of the integer programming problem $\eqref{IP}$ and the $l_0$-norm minimization problem $\eqref{la1}$ respectively.
\end{theorem}

\begin{proof}  For any optimal solution $x^*$ of problem $\eqref{la1}$, let $\hat{x}^*=(\lceil x_1^* \rceil, \lceil x_2^* \rceil, \cdots,$ $\lceil x_n^* \rceil)^T\in\{0,1\}^n$. By noting that $\hat{x}^*\geq x^*$ and $A\geq 0$, we have $A\hat{x}^*\geq Ax^* \geq b$. Hence, $\hat{x}^*$ is a feasible solution of problems $\eqref{IP}$ and $\eqref{la1}$ respectively. Further, since $\sum\limits_{i=1}^{n}|x^*_i|_0=\sum\limits_{i=1}^{n}|\lceil x_i^* \rceil|_0=\sum\limits_{i=1}^{n}|\hat{x}^*_i|_0$,  $\hat{x}^*$ is an optimal solution of problem $\eqref{la1}$.

Next, for any $x\in \{x\in \{0,1\}^n: Ax\ge b\}$, $x$ is also a feasible solution of problem $\eqref{la1}$. Since $\sum\limits_{i=1}^{n}x_i=\sum\limits_{i=1}^{n}|x_i|_0$, and by $\sum\limits_{i=1}^{n}|x^*_i|_0=\sum\limits_{i=1}^{n}|\hat{x}^*_i|_0$, we can obtain that $\sum\limits_{i=1}^{n}x_i \geq \sum\limits_{i=1}^{n}|\hat{x}^*_i|_0=\sum\limits_{i=1}^{n}\hat{x}^*_i$. Hence $\hat{x}^*$ is also  an optimal solution of problem $\eqref{IP}$. $\hfill\square$
\end{proof}

Theorem \ref{thm21} implies that the following corollary holds.

\begin{corollary}
If problem $\eqref{la1}$ has a unique integer optimal solution $\hat{x}^*$, then $\hat{x}^*$ is also an optimal solution of the integer programming problem $\eqref{IP}$.
\end{corollary}

Conversely, we next prove that an optimal solution of problem $\eqref{IP}$ is also an optimal solution of problem $\eqref{la1}$.

\begin{theorem}\label{thm22}
If $x$ is an optimal solution of problem $\eqref{IP}$, then $x$ is also an optimal solution of problem $\eqref{la1}$.
\end{theorem}
\begin{proof}
If $x$ is an optimal solution of problem $\eqref{IP}$, then $x$ is a feasible solution of problem $\eqref{la1}$, and satisfies that $\sum\limits_{i=1}^{n}x_i=\sum\limits_{i=1}^{n}|x_i|_0$. For any optimal solution $x^*$ of problem $\eqref{la1}$, by Theorem \ref{thm21}, $\hat{x}^*=(\lceil x_1^* \rceil ,\lceil x_2^* \rceil ,\cdots,\lceil x_n^* \rceil)^T$ is an optimal solution of problems $\eqref{IP}$ and $\eqref{la1}$. Hence,  $\sum\limits_{i=1}^{n}x_i=\sum\limits_{i=1}^{n}\hat{x}^*_i=\sum\limits_{i=1}^{n}|\lceil x_i^* \rceil|_0=\sum\limits_{i=1}^{n}|x^*_i|_0$. So $x$ is also an optimal solution of problem $\eqref{la1}$. $\hfill\square$
\end{proof}

By introducing  slack variables, letting
$A^{\prime}: =[A_1,A_2]\in R^{(m+n)\times(m+2n)}$, where $A_1=\big(\begin{array}{c}A \\  I\end{array}\big)\in R^{(m+n)\times n}$,
$A_2 = \big(\begin{array}{cc}-I&0\\0&I\end{array}\big)\in R^{(m+n)\times (m+n)}$,  $b^\prime=\big(\begin{array}{c} b \\ 1 \end{array}\big)
\in R^{m+n}$, then problem (\ref{la1}) can be rewritten as the form
\begin{equation}\label{la2}
\min\limits_{x,y}\{\|x\|_0: A_1x+A_2y=b^\prime,x\geq 0,y\geq0\}.
\end{equation}
Correspondingly, the weighted linear programming problem (\ref{wlp}) can be written in the form of the partially weighted linear programming problem
\begin{equation}\label{la3}
\min\limits_{x,y}\{c^Tx: A_1x+A_2y=b^\prime,x\geq 0,y\geq0\}.
\end{equation}

Then, to study the equivalence between the integer programming problem \eqref{IP} and the weighted linear programming problem \eqref{wlp}, we turn to derive conditions under which problems \eqref{la2} and \eqref{la3} have the same optimal solutions. In the sequel, we adapt the concept of $s$-goodness \cite{Juditsky2011} to problem  \eqref{la2}, such that problems \eqref{la2} and \eqref{la3} have the same unique optimal solution. Through this optimal solution, we can obtain the optimal solution of problem \eqref{IP}.

\section{ Nonnegative partial $s$-goodness}\label{chaper3}
Since problem \eqref{la2} is NP-hard, we are  interested in establishing some conditions, under which both of problems \eqref{la2} and \eqref{la3} have the same unique optimal solution.
Firstly, we give the following nonnegative partial s-goodness definition of matrix $A^{\prime}$ and weighted vector $c$, where $A^\prime=(A_1, A_2)$.  

\subsection{Definition of nonnegative partial $s$-goodness}\label{dnps}

\begin{definition}\label{dingyi1}
Let $A^{\prime}$ be a $(m+n)\times (m+2n)$ matrix and $s$ be an integer, $0\le s \leq n$. $0< c_i\leq 1$, $i=1, 2, ..., n$. We say that $(A^{\prime}, c)$ is nonnegative partially $s$-good   with respect to the columns of $A_1$, if for any pair of vectors $w^1\geq0\in R^n$, $w^2\geq0\in R^{m+n}$ with that $w^1\in R^n$ has at most $s$ nonzero elements, $(w^1, w^2)^T$ is the unique optimal solution to the optimization problem
\begin{equation}\label{EQ}
\min_{x, y}\big\{c^Tx:A_1x+A_2y=A_1w^1+A_2w^2,x\geq 0,y\geq0\big\}.
\end{equation}
\end{definition}

For the convenient of description, we say $(A^{\prime}, c)$ is nonnegative partially $s$-good to refer that $(A^{\prime}, c)$ is nonnegative partially $s$-good  with respect to the columns of $A_1$. Moreover, without loss of generality, for $s\in \{0,1,2,\cdots,n\}$, we say $w^1\in R^n$ with $\|w^1\|_0\leq s$ to mean that the nonzero entries of $w^1$ is no more than $s$. Meanwhile, in this paper, a vector is said to be $s$-sparse when this vector contains at most $s$ nonzero components.

It is obvious that, if the partially weighted linear problem \eqref{la3} has multiple optimal solutions, then $(A^{\prime},c)$ is not nonnegative partially $s$-good. However, we want to recover the optimal solution of problem \eqref{la2} from problem \eqref{la3}, in which $c$ is not fixed. So we  adjust the coefficient $c$ in problem \eqref{la3}, such that one of the optimal solutions is the unique optimal solution of problem \eqref{la3}.

By Definition \ref{dingyi1}, taking a step closer to our goal, we obtain the following results, which characterize the consistency of solutions to the problems \eqref{la2} and \eqref{la3}.

\begin{theorem}\label{thmknownS}
For any optimal solution $(w^1, w^2)^T$ of problem $(\ref{la2})$, where $w^1$ is an $s$-sparse vector, if $(A^\prime, c)$ is nonnegative partially $s$-good, then $(w^1,w^2)^T$ is the unique optimal solution to the partially weighted linear programming problem $(\ref{la3})$.
\end{theorem}
\begin{proof}

For any optimal solution $(w^1, w^2)^T$ of problem $(\ref{la2})$, where $w^1$ is an $s$-sparse vector, it holds that $A_1x+A_2y=b^\prime=A_1w^1+A_2w^2$ and $w^1\geq 0$, $w^2\geq 0$. That is, $(w^1,w^2)^T$ is a feasible solution of problem $(\ref{la3})$. If $(A^\prime, c)$ is nonnegative partially $s$-good, then according to Definition \ref{dingyi1}, $(w^1,w^2)^T$  with $\|w^1\|_0\leq s$ is a unique optimal solution of problem $(\ref{la3})$. $\hfill\square$
\end{proof}

Under the nonnegative partial $s$-goodness condition, Theorem \ref{thmknownS} shows that an optimal solution of problem $(\ref{la2})$ is a unique optimal solution of problem $(\ref{la3})$. In the above theorem, there is no requirement for the uniqueness of the solution of the $l_0$-norm minimization problem  $(\ref{la2})$, while only the uniqueness of the partially weighted linear problem $(\ref{la3})$ is required.

Next,  we give a stronger result that both problems $\eqref{la2}$ and $\eqref{la3}$ have a unique optimal solution.

\begin{theorem}\label{thm32}
Given an integer $0\leq s\leq n$, and let $(w^1,w^2)^T$ with $||w^1||_0\le s$ be a feasible solution to problem $(\ref{la3})$. Suppose $(A^\prime, c)$ is nonnegative partially $s$-good, then $(w^1,w^2)^T$ is both the unique optimal solution to the partially weighted linear problem $(\ref{la3})$ and the $l_0$-norm minimization problem $(\ref{la2})$.
\end{theorem}

\begin{proof}
Suppose $(A^\prime,c)$ is nonnegative partially $s$-good, and $(w^1,w^2)^T$ with $||w^1||_0\le s$ is a feasible solution to problem $(\ref{la3})$. Then by Definition \ref{dingyi1}, $(w^1,w^2)^T$ is the unique optimal solution to problem $(\ref{la3})$.

Next, we prove that $(w^1,w^2)^T$ is also the unique optimal solution to the $l_0-$norm minimization problem (\ref{la2}). Suppose $(x^1,y^1)^T$ is another solution to problem (\ref{la2}), and  $\|x^1\|_0\leq \|w^1\|_0\leq s$. Then $A_1x^1+A_2y^1=b^\prime=A_1w^1+A_2w^2$ and $x^1\geq 0$, $y^1\geq 0$. By Definition \ref{dingyi1},
$(x^1, y^1)^T$ is also a unique optimal solution to problem (\ref{la3}), and then we can get $(x^1,y^1)^T=(w^1,w^2)^T$. Hence  $(w^1, w^2)^T$ is the unique optimal solution of the $l_0-$norm minimization problem (\ref{la2}). $\hfill\square$
\end{proof}

According  to Section \ref{chaper22} and Theorems \ref{thmknownS} and \ref{thm32}, we can immediately get that an optimal solution of the integer programming problem \eqref{IP} can be recovered from problem $(\ref{la3})$, as in the following corollary.

\begin{corollary}
Suppose $(A^\prime, c)$ is nonnegative partially $s$-good, let $(w^1, w^2)^T$ with $||w^1||_0$ $\le s$ be an optimal solution to the partially weighted linear programming problem $(\ref{la3})$. Then $\lceil w^1 \rceil$ is an optimal solution of integer programming problem \eqref{IP}.
\end{corollary}

It seems not easy to completely characterize the nonnegative partial $s$-goodness of the constraint matrix $A^\prime$ and the coefficient $c$ of the objective function \eqref{wlp}. In the next subsection, we utilize two quantities to characterize nonnegative partial $s$-goodness.

\subsection{Two quantities of nonnegative partial $s$-goodness}\label{chaper31}
In this section, we introduce two quantities: $\gamma_{s,K}\big(A^\prime,c,\beta\big)$ and $\hat\gamma_{s,K}\big(A^\prime,c,\beta\big)$, where $K:=\{1,2,\cdots,n\}$ is the index set of $x$, i.e., the index set of the columns of matrix $A_1$.
In particular, for a vector $\theta\in R^{m+n}$, let $\|\cdot\|_*$ be the dual norm of $\|\cdot\|$ specified by $\|\theta\|_*=\max\limits_{d}\{d^T\theta:\|d\|\leq1\}$. In this paper, we consider the dual norm of $\|\cdot\|_1$.

\begin{definition}\label{def2}
Let $A^\prime\in R^{(m+n)\times (m+2n)}$, $s$ is an integer and $0\leq s\leq n$, $0< c_i\leq 1$, $i=1, 2, \dots,n$, $\beta\in [0,\infty]$.
We define $\gamma_{s,K}\big(A^\prime,c,\beta\big),\ \hat\gamma_{s,K}\big(A^\prime,c,\beta\big)$ as follows:

(1) $\gamma_{s,K}\big(A^\prime,c,\beta\big)$ is the infimum of $\gamma>0$ such that for every pair of vectors $z^1\in R^n,z^2\in R^{m+n}$, where $z^1\in R^n$ has $s$ nonzero entries, each  is equal to 1, there exists a vector $\theta\in R^{m+n}$ such that
\begin{eqnarray} \label{gamma1}
\| \theta\|_*\leq\beta,~ \big(A_1^{ T}\theta\big)_i\left\{\begin{array}{cc}=c_iz^1_i,& ~if~ z^1_i=1;\\
 \in [-\gamma,\gamma],&~if~ z^1_i=0,\end{array}\right.
 ~\text{and}~ \big(A_2^{ T}\theta\big)_i\left\{\begin{array}{cc}
=0,&~if~ z^2_i\neq 0;\\
\leq 0,&~if~ z^2_i=0.
\end{array}\right.
\end{eqnarray}

If for some pair of vectors $z^1\in R^n,z^2\in R^{m+n}$ as above,  there does no $\theta$ with $\| \theta\|_*\leq\beta$, such that $A_1^{ T}\theta$ coincides with $c\circ z^1$ on the support set of $z^1$, and $A_2^{ T}\theta$ coincides with $0$ on the support set of $z^2$, then we let $\gamma_{s,K}\big(A^\prime,c,\beta\big)=+\infty$.

(2) $\hat\gamma_{s,K}\big(A^\prime,c,\beta\big)$ is the infimum of $\gamma>0$ such that for every pair of vectors $z^1\in R^n,z^2\in R^{m+n}$, where $z^1\in R^n$ has $s$ nonzero entries, each is equal to 1, there exists a vector $\hat{\theta}\in R^{m+n}$ such that
\begin{eqnarray} \label{gamma2}
\|\hat{\theta}\|_*\leq\beta,~\|\big(A_1^{T}\hat{\theta}\big)-c\circ z^1\|_{\infty}\leq \gamma~
\text{and}~ \big(A_2^{ T}\hat{\theta}\big)_i\left\{\begin{array}{cc}
=0,&~if~ z^2_i\neq 0;\\
\leq 0,&~if~ z^2_i=0,
\end{array}\right.
\end{eqnarray}
\end{definition}
here $c\circ z^1$ means entry-wise product of the two vectors.

Furthermore, when $\beta=\infty$, we write $\gamma_{s,K}\big(A^\prime, c\big)$, $\hat{\gamma}_{s,K}\big(A^\prime, c\big)$ instead of $\gamma_{s,K}\big(A^\prime, c,$ $ \infty\big)$ and $\hat{\gamma}_{s,K}\big(A^\prime,c,\infty\big)$, respectively.


\begin{remark}\label{re1}
 Obviously, the set of values of the $\gamma$ is closed. Thus, if $\gamma_{s,K}\big(A^\prime,c,\beta)<+\infty$, then for every pair of vectors $z^1\in R^n,z^2\in R^{m+n}$, where $z^1\in R^n$ has $s$ nonzero entries, each is equal to 1, there exists a vector $\theta\in R^{m+n}$ such that
 \begin{equation}\label{B}
 \begin{array}{ll}
\| \theta\|_*\leq\beta, & \big(A_1^{ T}\theta\big)_i\left\{\begin{array}{cc}=c_iz^1_i, & ~if~ z^1_i=1;\\
 \in [-\gamma_{s,K}\big(A^\prime,c,\beta),\gamma_{s,K}\big(A^\prime,c,\beta)],&~if~ z^1_i=0,\end{array}\right.
 \\
 \\
&\big(A_2^{ T}\theta\big)_i\left\{\begin{array}{cc}
=0,&~if~ z^2_i\neq0;\\
\leq 0,&~if~ z^2_i=0.\end{array}\right.
 \end{array}
\end{equation}
Similarly, for every pair of vectors $z^1\in R^n$, $z^2\in R^{m+n}$, where $z^1\in R^n$ has $s$ nonzero entries, each is equal to 1, there exists a vector $\theta\in R^{m+n}$ such that
\begin{equation}\label{C}
\|\hat{\theta}\|_*\leq\beta, ~\|\big(A_1^{T}\hat{\theta}\big)-c\circ z^1\|_{\infty}\leq \hat{\gamma}_{s,K}\big(A^\prime,c,\beta\big)~\text{and}~ \big(A_2^{ T}\hat{\theta}\big)_i\left\{\begin{array}{cc}
=0,&~if~ z^2_i\neq 0;\\
\leq 0,&~if~ z^2_i=0.
\end{array}\right.
\end{equation}
\end{remark}

Before characterizing nonnegative partial $s$-goodness of  $(A^{\prime}, c)$ more specifically, we need to give some basic properties of $\gamma_{s,K}\big(A^\prime,c,\beta\big)$ and $\hat{\gamma}_{s,K}\big(A^\prime,c,\beta\big)$, such as convexity and monotonicity. 
Since nonnegative partial $s$-goodness of $(A^\prime,c)$ requires $\gamma_{s,K}(A^\prime,c)<\infty$, we assume it holds without loss of generality in the sequel.


\begin{lemma}\label{lam2}
$\gamma_{s,K}(A^\prime,c,\beta)$ and $\hat{\gamma}_{s,K}(A^\prime,c,\beta)$ are convex nonincreasing function of $\beta\in [0,+\infty]$.
\end{lemma}
\begin{proof}
Here, we only need to prove that the $\gamma_{s,K}\big(A^\prime,c,\beta)$ is a convex nonincreasing function with respect to $\beta\in [0,+\infty]$.  The property with respect to  $\hat{\gamma}_{s,K}(A^\prime,c,\beta)$ can be proved similarly.

Firstly, for the given $A^\prime, c$ and $s$, we demonstrate that the $\gamma_{s,K}\big(A^\prime, c, \beta)$ is a nonincreasing function of $\beta$. For any $ \beta_2\ge \beta_1$, according to the definition of $\gamma_{s,K}\big(A^\prime,c,\beta\big)$ and Remark \ref{re1},
for every pair of vectors $z^1\in R^n,z^2\in R^{m+n}$, where $z^1\in R^n$ has $s$ nonzero entries, each is equal to 1, there exists a vector $\theta\in R^{m+n}$ such that
 \begin{equation*}
 \begin{array}{ll}
\| \theta\|_*\leq\beta_1 ~,&\big(A_1^{T}\theta\big)_i\left\{\begin{array}{cc}=c_iz^1_i,& ~if~ z^1_i=1;\\
 \in [-\gamma_{s,K}\big(A^\prime,c,\beta_1),\gamma_{s,K}\big(A^\prime,c,\beta_1)],&~if~ z^1_i=0,\end{array}\right.
 \\
 \\
 &\big(A_2^{ T}\theta\big)_i\left\{\begin{array}{cc}
=0,&~if~ z^2_i\neq0;\\
\leq 0,&~if~ z^2_i=0.\end{array}\right.
 \end{array}
\end{equation*}
Since $\beta_2\ge\beta_1$, the $\theta$ in the above equation  also satisfies that
\begin{equation*}
 \begin{array}{ll}
\| \theta\|_*\leq \beta_2 ~,&\big(A_1^{T}\theta\big)_i\left\{\begin{array}{cc}=c_iz^1_i,& ~if~ z^1_i=1;\\
 \in [-\gamma_{s,K}\big(A^\prime,c,\beta_1),\gamma_{s,K}\big(A^\prime,c,\beta_1)],&~if~ z^1_i=0,\end{array}\right.
 \\
\\&\big(A_2^{ T}\theta\big)_i\left\{\begin{array}{cc}
=0,&~if~ z^2_i\neq0;\\
\leq 0,&~if~ z^2_i=0.\end{array}\right.
 \end{array}
\end{equation*}
Hence by the definition of $\gamma_{s,K}\big(A^\prime,c,\beta_2\big)$, $\gamma_{s,K}\big(A^\prime,c,\beta_1\big)\geq\gamma_{s,K}\big(A^\prime,c,\beta_2\big)$.


Next, we  prove that $\gamma_{s,K}(A^\prime,c,\beta)$ is a convex function of $\beta$. That is to say,  for any $\beta_1$, $\beta_2\in[0,+\infty]$, for any $\alpha\in[0,1]$, we need to prove that
\begin{equation}\label{eq9}
\gamma_{s,K}\big(A^\prime,c,\alpha\beta_1+(1-\alpha)\beta_2)\leq\alpha \gamma_{s,K}(A^\prime,c,\beta_1)+(1-\alpha)\gamma_{s,K}(A^\prime,c,\beta_2).
\end{equation}

Note that, the above inequality (\ref{eq9}) obviously holds if one of $\beta_1$ and $\beta_2$ is $+\infty$. Therefore, we only need to verify that for $\beta_1, \beta_2\in [0,+\infty)$, the inequality (\ref{eq9}) still holds. By the definition of $\gamma_{s,K}\big(A^\prime,c,\beta)$, it is easy to know that for every pair of vectors $z^1\in R^n,z^2\in R^{m+n}$, where $z^1\in R^n$ has $s$ nonzero entries, each is equal to 1, there exists a vector $\theta_\ell\in R^{m+n}$, $\ell\in\{1,2\}$ such that
 \begin{equation*}
 \begin{array}{ll}
\|\theta_{\ell}\|_*\leq\beta_{\ell} ~,&(A_1^{T}\theta_\ell)_{i}\left\{\begin{array}{cc}=c_iz^1_i,& ~if~ z^1_i=1;\\
 \in [-\gamma_{s,K}\big(A', c, \beta_{\ell}), \gamma_{s,K}\big(A', c, \beta_{\ell})], &~if~ z^1_i=0,\end{array}\right.
  \\
\\&\big(A_2^{ T}\theta_\ell\big)_i\left\{\begin{array}{cc}
=0,&~if~ z^2_i\neq0;\\
\leq 0,&~if~ z^2_i=0.\end{array}\right.
 \end{array}
\end{equation*}

Clearly, for any $\alpha\in[0,1]$, we can easily get $$\|\alpha\theta_1+(1-\alpha)\theta_2\|_*\leq\alpha\beta_1+(1-\alpha)\beta_2.$$
Moreover,
 \begin{equation*}
 \begin{array}{ll}
&[A_1^T(\alpha\theta_1+(1-\alpha)\theta_2)]_i\left\{\begin{array}{cc}=c_iz^1_i,& ~if~ z^1_i=1;\\
 \in [-k\varrho,\varrho],&~if~ z^1_i=0,\end{array}\right.
  \\
\\&[A_2^{T}(\alpha\theta_1+(1-\alpha)\theta_2)\big]_i\left\{\begin{array}{cc}
=0,&~if~ z^2_i\neq0;\\
\leq 0,&~if~ z^2_i=0,\end{array}\right.
 \end{array}
\end{equation*}
where $\varrho=\alpha\gamma_{s,K}(A^\prime, c, \beta_1)+(1-\alpha)\gamma_{s,K}(A^\prime, c,\beta_2)$. Hence, by the definition of $\gamma_{s,K}(\cdot)$, 
it holds that
$$\gamma_{s,K}\big(A^\prime, c,\alpha\beta_1+(1-\alpha)\beta_2)\leq\alpha \gamma_{s,K}(A^\prime, c,\beta_1)+(1-\alpha)\gamma_{s,K}(A^\prime, c,\beta_2). \qquad\square $$
\end{proof}

By Definition \ref{def2} and Lemma \ref{lam2}, the set of values of the $\gamma$ is closed and has a infimum $\gamma_{s,K}(A^\prime, c,\beta)$. Namely, for the given $A^\prime, c$ and $s$, if $\beta$  is large enough,  we can set $\gamma_{s,K}(A^\prime,c,\beta)=\gamma_{s,K}(A^\prime,c)$. In the same way, for the given $A^\prime, c$ and $s$, if $\beta$  is large enough, we can set $\hat{\gamma}_{s,K}(A^\prime,c,\beta)=\hat{\gamma}_{s,K}(A^\prime,c)$.

From the definitions of $\gamma_{s,K}(A^\prime,c,\beta)$ and $\hat{\gamma}_{s,K}(A^\prime,c,\beta)$, it is obvious that $s$ is another important parameter of $\gamma$. Next,  we  give a  property of $\gamma_{s,K}(A^\prime,c,\beta)$ and $\hat{\gamma}_{s,K}(A^\prime,c,\beta)$ with respect to $s$.


\begin{lemma}\label{lam1}
$\gamma_{s,K}(A^\prime,c,\beta)$ and $\hat{\gamma}_{s,K}(A^\prime,c,\beta)$ are  monotonically nondecreasing functions of  the parameter $s$.
\end{lemma}
\begin{proof}
Firstly, we prove that $\gamma_{s,K}(A^\prime,c,\beta)$ is a monotonically nondecreasing function of  the parameter $s$. Let $\gamma_{s,K}(A^\prime,c,\beta)<\infty$. According to the definition of $\gamma_{s,K}(A^\prime,c,\beta)$ and Remark \ref{re1}, for every pair of vectors $z^1\in R^n$, $z^2\in  R^{m+n}$, where $z^1\in R^n$ has $s$ nonzero entries, each is equal to 1, there exists a vector $\theta\in R^{m+n}$ such that
 \begin{equation}\label{gammas}
 \begin{array}{ll}
\| \theta\|_*\leq\beta ~,&\big(A_1^{T}\theta\big)_i\left\{\begin{array}{cc}=c_iz^1_i,& ~if~ z^1_i=1;\\
 \in [-\gamma_{s,K}\big(A^\prime,c,\beta),\gamma_{s,K}\big(A^\prime,c,\beta)],&~if~ z^1_i=0,\end{array}\right.
 \\
 \\
 &\big(A_2^{ T}\theta\big)_i\left\{\begin{array}{cc}
=0,&~if~ z^2_i\neq0;\\
\leq 0,&~if~ z^2_i=0.\end{array}\right.
 \end{array}
\end{equation}
Then, let $t_s(z^1, z^2)$ be the minimal value of the optimization problem
\begin{equation}
\begin{array}{cl}\label{ops}
\min\limits_{\theta} & \|(A_1^T\theta)\|_\infty \\
s.t. & \|\theta\|_*\leq\beta \\
     & (A_1^T\theta)_i=c_iz^1_i, ~i\in I_1\\
     & (A_2^{T}\theta)_i=0,~ i\in I_2 \\
     & (A_2^{T}\theta)_i\leq 0, ~i\in \bar{I}_2,
\end{array}
\end{equation}
where $I_1=\{i: ~z^1_i=1\}$, $I_2=\{i: ~z^2_i\neq0\}$ and $\bar{I}_2=\{i:~z^2_i=0\}$. Obviously, $t_s(z^1, z^2)\le \gamma_{s,K}(A^\prime,c,\beta)$.


Let $s^\prime=s-1<s$. For every pair of vectors $z^{1\prime}\in R^n$, $z^2\in R^{m+n}$, where $z^{1\prime}\in R^n$ has $s-1$ nonzero entries, each is equal to 1, we can construct a pair of vectors $z^1\in R^n$, $z^2\in  R^{m+n}$, where $z^1$ is obtained from $z^{1\prime}$ by changing one entry with value $0$ to $1$. According to \eqref{ops}, it is obvious that
    $$t_{s'}(z^{1\prime}, z^2)\le t_s(z^1, z^2)\le \gamma_{s,K}(A^\prime,c,\beta). $$
Hence, $\gamma_{s', K}(A^\prime, c, \beta)\le \gamma_{s, K}(A^\prime,c,\beta)$.

Using the way similar to the above proof,
we can also show that $\hat{\gamma}_{s^\prime,K}(A^\prime,c,\beta)\leq \hat{\gamma}_{s,K}(A^\prime,c,\beta)$, i.e.,  $\hat{\gamma}_{s,K}(A^\prime, c, \beta)$ is a monotonically nondecreasing function of $s$. $\hfill\square$

\end{proof}
\begin{remark} According to Lemma \ref{lam1},
$\gamma_{s,K}(A^\prime,c,\beta)$ and $\hat{\gamma}_{s,K}(A^\prime,c,\beta)$ are  nondecreasing functions of $s$, hence for all $s'\le s$, Remark \ref{re1} holds. 
That means, in Remark \ref{re1}  for every pair of vectors $z^1\in R^n$, $z^2\in  R^{m+n}$ with $\|z^1\|_0\le s$, there exists $\theta$ with $\|\theta\|_*\leq\beta$ such that Equations  \eqref{gamma1} and  \eqref{gamma2} hold.
\end{remark}

\subsection{Sufficient condition and necessary condition of nonnegative partial $s$-goodness}\label{chaper32}
In this subsection, via $\gamma_{s,K}(A^\prime,c)$ we propose a sufficient condition and a necessary condition for nonnegative partial $s$-goodness of the constraint matrix $A^\prime$ and the coefficient $c$ of the objective function.
\begin{theorem}\label{EPI1}
Given $A^\prime\in R^{(m+n)\times (m+2n)}$, $s$ is an integer and $0\leq s\leq n$, $0< c\leq 1$, we can obtain:

 $(a)$ if $(A^\prime,c)$ is nonnegative partially $s$-good, then $\gamma_{s,K}(A^\prime,c)\leq\max\limits_{0< i\leq n}c_i$;

 $(b)$ if $\gamma_{s,K}(A^\prime,c)<\min\limits_{0< i\leq n} c_i$, then $(A^\prime,c)$ is nonnegative partially $s$-good.
\end{theorem}
\proof
$(a)$ Suppose $(A^\prime, c)$ is nonnegative partially $s$-good.
For any given $w=(w^1, w^2)^T$\\$\geq 0\in R^{n}\times R^{m+n}$ with $\|w^1\|_0\leq s$,
let $I_1=\{i:~w^1_i>0\}$, $\bar{I}_1=\{i:~w^1_i=0\}$, $I_2=\{i:~w^2_i>0\}$ and $\bar{I}_2=\{i:~w^2_i=0\}$.

By Definition \ref{dingyi1}, $w$ is the unique optimal solution to  problem ($\ref{EQ}$). According to the optimality condition, there exists $\theta\in R^{m+n}$ such that $f_\theta(x, y)=\sum\limits_{i=1}^{n}c_i|x_i|-\theta^T(A_1x+A_2y-A_1w^1-A_2w^2)$ attains its minimum value at $(x,y)^T=(w^1,w^2)^T$, i.e., $0\in\partial f_\theta(w^1, w^2)$. This implies  that 
 \begin{eqnarray*}
 \begin{array}{ll}
(A_1^T\theta)_i\left\{\begin{array}{cc}=c_i,& ~ i\in I_1;\\
\in[-\max\limits_{0< i\leq n}c_i, \max\limits_{0< i\leq n}c_i],& ~i\in \bar{I}_1,\end{array}\right.
~and~(A_2^{T}\theta)_i\left\{\begin{array}{cc}
=0,& ~i\in I_2;\\
\leq 0,&~i\in \bar{I}_2.\end{array}\right.
 \end{array}
\end{eqnarray*}

Since $w\geq 0$, hence, for the optimization problem
\begin{eqnarray*}
\min\limits_{\theta,
\gamma}\Bigg\{\gamma:(A_1^T\theta)_i\left\{\begin{array}{cc}=c_i,&~i\in I_1;\\
 \in [-\gamma,\gamma],&~i\in \bar{I}_1,\end{array}\right.
 ~and~(A_2^{T}\theta)_i\left\{\begin{array}{cc}
=0,&~ i\in I_2;\\
\leq 0,&~i\in \bar{I}_2,\end{array}\right. \Bigg \}
\end{eqnarray*}
the optimal value   $\gamma\leq\max\limits_{0< i\leq n}c_i$. By Definition \ref{def2}, $\gamma_{s,K}(A^\prime, c)$ is the infimum of $\gamma$, thus $\gamma_{s,K}(A^\prime,c)\leq\max\limits_{0< i\leq n}c_i$.

$(b)$ Suppose $\gamma_{s,K}(A^\prime,c)<\min\limits_{0< i\leq n} c_i$, next we  prove $(A^\prime,c)$ is nonnegative partially $s$-good. That is, for a vector $w=(w^1,w^2)^T\geq0$ with $A_1w^1+A_2w^2=b^\prime$ and $\|w^1\|_0\leq s$, we need to prove that $w$ is the unique optimal solution to  problem ($\ref{EQ}$).

First, we consider the special case that $(w^1,w^2)^T=(0, w^2)^T$. Obviously, $(x,y)^T=(0,w^2)^T$ is the unique optimal solution to  problem ($\ref{EQ}$). The reason is that, when $w^1=0$, then $A_2y=A_2w^2$, and since $A_2 = \big(\begin{array}{ll}-I&0\\0&I\end{array}\big)$ is injective, it is easy to see that $y=w^2$ is unique.

Now, suppose $\|w^1\|_0=s^\prime$, $0\leq s^\prime\le s$, and its support set is $I_1$. Meanwhile, let $I_2$ be the support index set of $w^2$.
According to Lemma \ref{lam1}, we have $\gamma:=\gamma_{s^{\prime},K}(A^\prime,c)\leq\gamma_{s,K}(A^\prime,c)$. Since $\gamma_{s,K}(A^\prime,c)<\min\limits_{0< i\leq n}c_i$,  we get that $\gamma <\min\limits_{0< i\leq n}c_i$. Moreover, by the definition of $\gamma_{s,K}(\cdot)$, there exists $\theta\in R^{m+n}$ such that
 \begin{eqnarray}\label{eq16}
 \begin{array}{ll}
 \|\theta\|_*\leq\beta,~
(A_1^T\theta)_i\left\{\begin{array}{cc}=c_isign(w^1_i),& ~ i\in I_1;\\
\in [-\gamma,\gamma],& ~i\in \bar{I}_1,\end{array}\right.
~(A_2^{T}\theta)_i\left\{\begin{array}{cc}
=0,& ~i\in I_2;\\
\leq 0,&~i\in \bar{I}_2,\end{array}\right.
 \end{array}
\end{eqnarray}
where $I_1=\{i:~w^1_i>0\}$, $\bar{I}_1=\{i:~w^1_i=0\}$, $I_2=\{i:~w^2_i>0\}$ and $\bar{I}_2=\{i:~w^2_i=0\}$. Furthermore by (a),  there exists $\theta\in R^{m+n}$ satisfying Eq. \eqref{eq16} which is the optimal Lagrange multiplier of problem (\ref{EQ}). Then for any feasible solution $(x, y)^T$ of problem (\ref{EQ}), it holds that
\begin{equation}\label{f}
\begin{aligned}
f(x,y)&=c^Tx-\theta^T(A_1x+A_2y-A_1w^1-A_2w^2)\\
&=c^Tx-(A_1^T\theta)^T(x-w^1)-(A_2^T\theta)^T(y-w^2)\\
&=\sum_{i\in I_1}c_iw^1_i+\sum_{i\in \bar{I}_1}(c_i-(A_1^T\theta)_i)x_i-\sum_{i\in \bar{I}_2}(A_2^{T}\theta)_iy_i\\
&\geq\sum_{i\in I_1}c_iw^1_i+\sum_{i\in \bar{I}_1}(c_i-(A_1^T\theta)_i)x_i\\
&\geq c^Tw^1.
\end{aligned}
\end{equation}

According to \eqref{f}, it is obvious that the minimum value of the Lagrange function can be attained at $x=w^1$. Further, since $(x,y)$ and $(w^1,w^2)$ have the relationship  $A_1x+A_2y=A_1w^1+A_2w^2$,  $A_1x=A_1w^1$ and $A_2=\big(\begin{array}{ll}-I&0\\0&I\end{array}\big)$ is an injective matrix, it is easy to show that $A_2y=A_2w^2$. Therefore, $(x,y)^T=(w^1,w^2)^T$ is an optimal solution of problem (\ref{EQ}).

Next, it is necessary to prove that this optimal solution is unique. Suppose $(\tilde{x},\tilde{y})^T$ is another optimal solution of  problem (\ref{EQ}). That is, \begin{equation*}
\begin{aligned}
f(\tilde{x},\tilde{y})-f(w^1,w^2)=\sum\limits_{i\in \bar{I}_1}(c_i-(A_1^T\theta)_i)\tilde{x}_i-\sum\limits_{i\in \bar{I}_2}(A_2^{T}\theta)_i\tilde{y}_i=0.
\end{aligned}
\end{equation*}
By the assumption that $\gamma<\min\limits_{0< i\leq n}c_i$, and
from Eq. \eqref{eq16}, $|A_1^T\theta|_{i}< \min\limits_{0< i\leq n} c_i$ for all $i\in \bar {I}_1$, which means that   $\tilde{x}_i=0$ and $\sum\limits_{i\in \bar{I}_2}(A_2^{T}\theta)_i\tilde{y}_i=0$.
Therefore, $\tilde{x}_i=w^1_i=0$ for all $i\in \bar{I}_1$. Hence, $\|\tilde{x}-w^1\|_0\leq s$.

Further, for the above vector $\tilde{x}_i-w^1_i$, define
\begin{equation*}
\begin{aligned}
 h(\tilde{x}-w^1, \tilde{y}-w^2):=\sum\limits_{i=1}^{n}c_i|(\tilde{x}-w^1)_i|-\tilde\theta^T(A_1(\tilde{x}-w^1)+A_2(\tilde{y}-w^2)).
\end{aligned}
\end{equation*}
Similar to the proof of part (a) in Theorem \ref{EPI1}, there exists $\tilde\theta\in R^{m+n}$ such that 
 \begin{eqnarray*}
 \begin{array}{ll}
&(A_1^T\tilde\theta)_i\left\{\begin{array}{cc}=c_isign((\tilde{x}-w^1)_i),& ~ if ~ (\tilde{x}-w^1)_i\neq 0;\\
\in [-\max\limits_{0< i\leq n}c_i, \max\limits_{0< i\leq n}c_i],& ~if~ (\tilde{x}-w^1)_i= 0,\end{array}\right.
 \end{array}
\end{eqnarray*}
and
 \begin{eqnarray*}
 \begin{array}{ll}
(A_2^{T}\tilde\theta)_i\left\{\begin{array}{cc}
=0,& ~if ~(\tilde{y}-w^2)_i\neq 0;\\
\leq 0,&~if~ (\tilde{y}-w^2)_i= 0.\end{array}\right.
 \end{array}
\end{eqnarray*}

Therefore, for the $\tilde{\theta}$ in the function $h(\tilde{x}-w^1, \tilde{y}-w^2)$,  we have
\begin{equation*}
\begin{aligned}
0&=(A_1^T\tilde{\theta})^T(\tilde{x}-w^1)+(A_2^T\tilde{\theta})^T(\tilde{y}-w^2)\\
&=\sum\limits_{i\in I_1}(A_1^T\tilde{\theta})^T_i(\tilde{x}_i-w^1_i)+\sum\limits_{i\in I_2}(A_2^T\tilde{\theta})^T_i(\tilde{y}_i-w^2_i),
\end{aligned}
\end{equation*}
and then  we can get $\tilde{x}_i=w^1_i$ for all $i\in I_1$. This  combined with the fact $A_1\tilde{x}+A_2\tilde{y} =A_1w^1+A_2w^2$ can lead to $A_2\tilde{y} = A_2w^2$. Further, since  $A_2=\big(\begin{array}{ll}-I&0\\0&I\end{array}\big)$ is an injective matrix, we have $\tilde{y} =w^2$ and then $(\tilde{x},\tilde{y})^T=(w^1,w^2)^T$.$\hfill\square$


Below, we show   the relationship between  $\gamma_{s,K}(\cdot)$ and $\hat\gamma_{s,K}(\cdot)$.

\begin{proposition}\label{pro1}
For arbitrary $\beta\in[0,\infty]$, if $\hat{\gamma}:=\hat{\gamma}_{s,K}(A^\prime,c,\beta)<\frac{1}{2}\min\limits_{0< i\leq n}c_i$, then
\begin{eqnarray}\label{eqpro1}
\gamma_{s,K}(A^\prime,c,\frac{\min\limits_{0<i\leq n}c_i}{\min\limits_{0<i\leq n}c_i-\hat\gamma}\beta)\leq\frac{\min\limits_{0<i\leq n}c_i}{\min\limits_{0<i\leq n}c_i-\hat\gamma}\hat\gamma<\min\limits_{0<i\leq n}c_i.
\end{eqnarray}
\end{proposition}
\begin{proof}

Suppose $\hat\gamma:=\hat\gamma_{s,K}(A^\prime,c,\beta)<\frac{1}{2}\min\limits_{0< i\leq n}c_i$. Now let $I_1$ be an $s$-element subset of $\{1, 2, \dots, n\}$, $\bar{I}_1:=\{1, 2, \dots, n\}\setminus I_1$, and let $I_2$ be a subset of $\{n+1, n+2, \dots, m+2n\}$, $\bar{I}_2:=\{n+1, n+2, \dots, m+2n\}\setminus I_2$.

For the $I_1$, $\bar{I}_1$ and $I_2$, $\bar{I}_2$, we define a closed convex set $\Pi_{I_1}$ in $R^n$ as
\begin{small}
\begin{eqnarray*}
\begin{array}{ll}
\Pi_{I_1}=\left\{\tau^\prime\in R^{n}:\exists \theta\in R^{m+n},\|\theta\|_*\leq\beta,(A_1^{T}\theta)_i\left\{\begin{array}{cc}=c_i\tau^\prime_i,~i\in I_1;\\ \in [-\hat\gamma,
\hat\gamma],~i\in\bar{I}_1,\end{array}\right.
(A_2^{ T}\theta)_i\left\{\begin{array}{cc}
=0,&~i\in I_2;\\
\leq 0,&~i\in \bar{I_2}
\end{array}\right.\right\}.
\end{array}
\end{eqnarray*}
\end{small}

Similar to the proof of Proposition 2.1 in \cite{Juditsky2011}, we claim that $\Pi_{I_1}$ contains the $\|\cdot\|_\infty$-ball $B$ centered at the origin and with the radius $\frac{\min\limits_{0<i\leq n}c_i-\hat{\gamma}}{\min\limits_{0<i\leq n}c_i}$. The proof is as follows.

Define a subspaces $L_{I_1}:=\{\tau^\prime\in R^n:\tau^\prime_i=0,i\in \bar{I}_{1}\}$ and let $L^{\perp}_{I_1}:=\{\tau^\prime\in R^n:\tau^\prime_i=0,i\in I_{1}\}$ be the orthogonal
complement of $L_{I_1}$. Let $P$ be the projection of $\Pi_{I_1}$ onto $L_{I_1}$ and $P^\prime$ be the projection of $\Pi_{I_1}$ onto $L^{\perp}_{I_1}$. Note that $\Pi_{I_1}$ is the direct sum of $P$ and $P^\prime$. Thus, $P$ is a closed convex set. Obviously, $L_{I_1}$ can be naturally identified with $R^s$.

Hence, the claim in the above can be more precisely described as that the image $\bar{P}\subset R^s$ of $P$ contains the $\|\cdot\|_\infty$-ball $B_s$ centered at the origin and with the radius $\frac{\min\limits_{0<i\leq n}c_i-\hat{\gamma}}{\min\limits_{0<i\leq n}c_i}$ in $R^s$.

Next, we  prove that $B_s\subseteq \bar{P}$. Contradictorily, suppose that $\bar{P}$ does not contain $B_s$. Since $P$ is a closed convex set,  $\bar{P}$ is also a closed convex set.   According to the separating hyperplane theorem, for $\nu\in B_s\setminus \bar{P}$ and $\bar{\nu}\in \bar{P}$, there exists $u\in R^s$ with $\|u\|_1=1$, such that

\begin{eqnarray}\label{SHT}
\begin{array}{ll}
 u^T\nu>\max\limits_{\bar{\nu}\in \bar{P}}~ u^T\bar{\nu}.
\end{array}
\end{eqnarray}
In the following, we prove that there does not exist $u$ such that the inequality \eqref{SHT} holds.

First, define $\bar{z}^1\in R^n$ with $\|\bar{z}^1\|_0=s$ and $\bar{z}^2\in R^{m+n}$ as
 \begin{eqnarray*}
  \begin{array}{ll}
 \bar{z}^1=\left\{\begin{array}{cc}1,~i\in I_1;\\ 0,~i\in\bar{I}_1,\end{array}\right.
 ~\bar{z}^2\left\{\begin{array}{cc}\neq0,~i\in I_2;\\ =0,~i\in\bar{I}_2.\end{array}\right.
  \end{array}
 \end{eqnarray*}
By the definition of $\hat{\gamma}_{s,K}(A^\prime, c, \beta)$, for the given $\bar{z}^1\in R^n$, $\bar{z}^2\in R^{m+n}$, there exists $\bar{\theta}\in R^{m+n}$ such that

\begin{eqnarray*}
\begin{array}{ll}
\|\bar{\theta}\|_*\leq\beta,~\|\big(A_1^{T}\bar{\theta})-c^T\circ\bar{z}^1\|_\infty\leq \hat\gamma,~ \big(A_2^{ T}\theta\big)_i\left\{\begin{array}{cc}
=0,&~if~ \bar{z}^2_i\neq 0;\\
\leq 0,&~if~ \bar{z}^2_i=0.
\end{array}\right.
 \end{array}
 \end{eqnarray*}
Then for the $u$ with $\|u\|_1=1$ in  \eqref{SHT}, combining the above inequality with the definitions of $\Pi_{I_1}$ and $\bar{P}$,  there exists a vector $\nu^\prime\in\bar{P}$ such that $$|c_i\nu^\prime_i-c_isign(u_i)|\leq\hat{\gamma}, ~i\in I_1.$$
Thus,
$$sign(u_i)-\frac{\hat{\gamma}}{\min\limits_{0<i\leq n}c_i}\leq\nu^\prime_i\leq sign(u_i)+\frac{\hat{\gamma}}{\min\limits_{0<i\leq n}c_i}, ~i\in I_1.$$ Since $\hat{\gamma}<\frac{1}{2}\min\limits_{0<i\leq n}c_i$, the above inequalities  imply that the sign of $\nu_i^\prime$ is same as $u_i$, for all $i$. Moreover, according to the definition of $\hat{\gamma}_{s,K}(A^\prime, c, \beta)$, we can get $$1-\frac{\hat{\gamma}}{\min\limits_{0<i\leq n}c_i}\leq\nu^\prime_i\leq 1+\frac{\hat{\gamma}}{\min\limits_{0<i\leq n}c_i}, ~i\in I_1.$$
Hence, $\nu^\prime>0$,  $u>0$, and
$$\nu^\prime_i\geq\frac{\min\limits_{0<i\leq n}c_i-\hat{\gamma}}{\min\limits_{0<i\leq n}c_i}, i\in I_1.$$
So
\begin{eqnarray*}
u^T\nu^\prime\geq\sum\limits_{i=1}^{s}|u_i|\frac{\min\limits_{0<i\leq n}c_i-\hat{\gamma}}{\min\limits_{0<i\leq n}c_i}=\frac{\min\limits_{0<i\leq n}c_i-\hat{\gamma}}{\min\limits_{0<i\leq n}c_i}.
 \end{eqnarray*}
Further, for the above given $\nu\in B_s$ and $\|u\|_1=1$, we have
\begin{eqnarray*}
\frac{\min\limits_{0<i\leq n}c_i-\hat{\gamma}}{\min\limits_{0<i\leq n}c_i}\geq\|\nu\|_\infty=\|u\|_1\|\nu\|_\infty\geq u^T\nu>u^T\nu^\prime\geq\frac{\min\limits_{0<i\leq n}c_i-\hat{\gamma}}{\min\limits_{0<i\leq n}c_i},
\end{eqnarray*}
which is a contradiction. So the claim holds.

Through the above proof, we can conclude that, for any $z=(z^1,z^2)^T\in R^n\times R^{m+n}$ with $z^1_i=1$ for $i\in I_1$, and  $z^1_i=0$ else, there exists $\tau^\prime\in\Pi_{I_1}$ such that $$\tau^\prime_i=(\frac{\min\limits_{0<i\leq n}c_i-\hat{\gamma}}{\min\limits_{0<i\leq n}c_i}) z_i^1, ~i\in I_1.$$

Further, by the definition of $\Pi_{I_1}$, there exists $\hat\theta\in R^{m+n}$ with $\|\hat\theta\|_*\leq \frac{\min\limits_{0<i\leq n}c_i}{\min\limits_{0<i\leq n}c_i-\hat{\gamma}}\beta$ such that
 \begin{eqnarray*}
 \begin{array}{ll}
&(A_1^T\hat\theta)_i\left\{\begin{array}{cc}=\frac{\min\limits_{0<i\leq n}c_i}{\min\limits_{0<i\leq n}c_i-\hat\gamma}c_i\tau^\prime_i=c_iz^1_i,& ~ i\in I_1;\\
\in [-\frac{\min\limits_{0<i\leq n}c_i}{\min\limits_{0<i\leq n}c_i-\hat\gamma}\hat\gamma,\frac{\min\limits_{0<i\leq n}c_i}{\min\limits_{0<i\leq n}c_i-\hat\gamma}\hat\gamma],& ~ i\in\bar{ I}_1,\end{array}\right.\\
 \end{array}
\end{eqnarray*}
and
 \begin{eqnarray*}
 \begin{array}{ll}
&(A_2^{T}\hat\theta)_i\left\{\begin{array}{cc}
=0,& ~i\in I_2;\\
\leq 0,&~ i\in\bar{I}_2.\end{array}\right.
 \end{array}
\end{eqnarray*}
So by the definition of $\gamma_{s,K}(A^\prime,c,\beta)$, and since $\hat{\gamma}<\frac{1}{2}\min\limits_{0<i\leq n}c_i$, we can obtain

\begin{eqnarray*}
\gamma_{s,K}(A^\prime,c,\frac{\min\limits_{0<i\leq n}c_i}{\min\limits_{0<i\leq n}c_i-\hat\gamma}\beta)\leq\frac{\min\limits_{0<i\leq n}c_i}{\min\limits_{0<i\leq n}c_i-\hat\gamma}\hat\gamma<\min\limits_{0<i\leq n}c_i. \qquad \hfill\square
\end{eqnarray*} 
\end{proof}

Based on Proposition \ref{pro1}, Theorem \ref{EPI1} can be equivalently written as: 
\begin{theorem}\label{thm8}
Given $A^\prime\in R^{(m+n)\times (m+2n)}$, $s$ is an integer and $0\leq s\leq n$,  $0< c\leq 1$, if $\hat{\gamma}_{s,K}(A^\prime,c)<\frac{1}{2}\min\limits_{0< i\leq n} c_i$, then $(A^\prime,c)$ is nonnegative partially $s$-good.
\end{theorem}

According to Theorem \ref{EPI1}, to show that $(A^\prime, c)$ is nonnegative partially $s$-good, we need to compare the magnitude of $\gamma_{s,K}(A^\prime, c)$ with $\min\limits_{0< i\leq n}c_i$.   
Now, due to Theorem \ref{thm8}, and we know that $\hat{\gamma}_{s,K}(A^\prime, c, \beta)$ is weaker than $\gamma_{s,K}(A^\prime, c, \beta)$, hence we focus on $\hat\gamma(\cdot)$, which has a specific representation presented in the next subsection.


\subsection{Specific representation of $\hat{\gamma}_{s,K}(\cdot)$}\label{chaper33}

$\hat{\gamma}_{s,K}(A^\prime, c, \beta )$is given in Definition \ref{def2}, which is essentially obtained from the optimality condition of problem (\ref{EQ}).
In this subsection, we give a specific representation of $\hat{\gamma}_{s,K}(A^\prime, c, \beta)$ in more detail.


\begin{theorem}\label{EPI2}
Consider the polytope
\begin{equation*}
P_s=\{\tau\in R_+^{m+2n}: \tau=(\tau^1,\tau^2)^T,  \tau^1\in R^n,  \tau^2\in R^{m+n},~\|\tau^1\|_1\leq s,~\|\tau^1\|_\infty \leq 1\},
\end{equation*}
we have
\begin{equation}\label{F}
\hat{\gamma}_{s,K}(A^\prime,c,\beta)=\max_{\tau, x}\{\sum_{i=1}^{n} \tau^1_ic_ix_i-\beta\|A_1x\|_1:\tau\in P_s,\|x\|_1\leq1,~x\geq0\}.
\end{equation}
Particularly,
\begin{equation}\label{G}
\hat{\gamma}_{s,K}(A^\prime,c)=\max\limits_{\tau, x}\{\sum_{i=1}^{n} \tau^1_ic_ix_i:~\tau\in P_s,~\|x\|_1\leq 1,~A_1x=0,~x\geq0\}.
\end{equation}
\end{theorem}

\begin{proof}

For any vector $y\ge 0$, let $I_2(y)$ be its support set. According to Definition \ref{def2}, define
 $$B_\beta(y)=\{\theta\in R^{m+n}:~\|\theta\|_*\leq\beta,~(A_2^{T}\theta)_i=0 \text{ for } i\in I_2(y),~(A_2^{T}\theta)_i\leq0 \text{ otherwise}\},$$
 $$B=\{\nu\in R^{n}:~\|\nu\|_\infty\leq 1\}.$$
 Then, $\hat{\gamma}_{s,K}(A^\prime,c,\beta)$ is the smallest $\gamma$, such that the set $C_{1,\gamma,\beta}:=A_1^{T}B_\beta(y)+\gamma B\subseteq R^n$ is closed,  convex,  and contains all vectors with $s$ nonzero elements, which are selected from $c_i$, $i=1, 2, \cdots, n$. This statement is equivalent to that $C_{1,\gamma,\beta}$ contains the convex hull of the vectors.

Let $\hat{c}=(c^T, 0, \cdots, 0)^T\in R^{m+2n}$. Then $C_{1,\gamma,\beta}$ contains the projection of $\hat{c}\circ P_s$ onto the $R^n$ space. Thus, $\gamma$ satisfies the relationship of $\hat{c}\circ P_s\subseteq C_{1,\gamma,\beta}\times R^{m+n}$, if and only if for any pair of $(x,y)^T\geq0$ with $x\in R^n$ and $y\in R^{m+n}$,
\begin{equation}\label{UP}
\begin{aligned}
&\max_{\tau\in P_s}\sum_{i=1}^{n} \tau^1_ic_ix_i \leq\max_{\eta\in C_{1,\gamma,\beta}}\{\sum_{i=1}^{n} \eta_ix_i\}\\
&=\max_{\theta,\nu}\{<x,A_1^{T}\theta>+\gamma <x,\nu>: ~\theta\in B_\beta(y),~ \|\nu\|_\infty\leq 1
\}\\
&\le\max_{\theta,\nu}\{<x,A_1^{T}\theta>+\gamma <x,\nu>: ~\|\theta\|_*\leq \beta,~ \|\nu\|_\infty\leq 1\}\\
&=\max_{\theta,\nu}\{<A_1x,\theta>+\gamma <x,\nu>: ~\|\theta\|_*\leq \beta,~ \|\nu\|_\infty\leq 1\}\\
&=\beta\|A_1x\|_1+\gamma\|x\|_1.
\end{aligned}
\end{equation}
That is, $\hat{c}\circ P_s\subseteq C_{1,\gamma,\beta}\times R^{m+n}$ if and only if $\gamma$ satisfies that
\begin{equation*}
\begin{aligned}
\max_{\tau\in P_s}~ \{\sum_{i=1}^{n} \tau^1_ic_ix_i-\beta\|A_1x\|_1:~x\geq0\}\leq \gamma\|x\|_1,
\end{aligned}
\end{equation*}
namely,
\begin{equation*}
\begin{aligned}
\max_{\tau, x}\{\sum_{i=1}^{n} \tau^1_ic_ix_i-\beta\|A_1x\|_1:\tau\in P_s,\|x\|_1\leq1,~x\geq0\}\leq\gamma.
\end{aligned}
\end{equation*}
Therefore Eq. \eqref{F} holds, since $\hat{\gamma}_{s,K}(A^\prime,c,\beta)$ is the smallest $\gamma$. Finally, Eq. \eqref{G} holds due to that Eq. \eqref{F} can be regarded as a penalty problem of  Eq. \eqref{G}.  $\hfill\square$
\end{proof}


For $c\circ x\geq0\in R^{n}$, we define the sum of the $s$ largest entries of $c\circ x$ as $$\|c\circ x\|_{s,K,1}:= \max\limits_{\tau\in P_s}\sum\limits_{i=1}^{n} \tau^1_ic_ix_i.$$
Then by taking Theorems \ref{thm8} and \ref{EPI2} into consideration, we can obtain the following result:
\begin{corollary}\label{CO1}
Given a matrix $A_1$, $\hat\gamma_{s,K}(A^\prime,c)$ is the least upper bound on $\|c\circ x\|_{s,K,1}:=\max\limits_{\tau\in P_s}\sum\limits_{i=1}^{n} \tau^1_ic_ix_i$ over $x\geq0$, $x\in Ker(A_1)$ and $\|x\|_1\leq1$. As a result, if the maximum of $\|c\circ x\|_{s,K,1}$ over $x\in Ker(A_1)$ and $\|x\|_1\leq1$ is less than $\frac{1}{2}\min\limits_{0<i\leq n}c_i$, then $(A^\prime,c)$ is nonnegative partially $s$-good.
\end{corollary}

Equations \eqref{F} and \eqref{G} provide the specific forms of  $\hat\gamma_{s, K}(\cdot)$. Thus,  we can judge the nonnegative partial $s$-goodness of $(A^\prime,c)$ according to Theorem \ref{thm8}, as long as the Eqs. \eqref{F} and \eqref{G} can be calculated. However, in Eqs. \eqref{F} and \eqref{G}, the calculation of $\hat\gamma_{s,K}(\cdot)$ is complicated, and sometimes it is not easy to directly calculate the specific value of $\hat\gamma_{s,K}(A^\prime,c,\beta)$.
In order to make up for this shortcoming, in what follows, we give  effective upper bounding of $\hat\gamma_{s,K}(A^\prime,c,\beta)$ to estimate the value of $\hat\gamma_{s,K}(A ^\prime,c,\beta)$.


\section{ Efficient bounding of $\hat\gamma_{s,K}(\cdot)$}\label{chaper4}
From the previous sections, we have shown that $\hat\gamma_{s,K}(A^\prime,c,\beta)$ plays an important role  in distinguishing whether $(A^\prime, c)$ is non-negative partially $s$-good. However, it still not easy to determine the exact value of $\hat\gamma_{s,K}(A^\prime,c,\beta)$ according to Eqs. \eqref{G} and \eqref{F}. In this section, we  introduce  efficiently computable upper bound on the value of $\hat\gamma_{s,K}(A^\prime,c,\beta)$.

Since $\hat\gamma_{s,K}(A^\prime,c,\beta)\geq\hat\gamma_{s,K}(A^\prime,c)$ for any $\beta>0$, we will use Eq. $\eqref{G}$ to calculate an upper bound of $\hat\gamma_{s,K}(A^\prime,c)$. The difficulty of calculation is in the process of the linear constraint $A_1x=0$, which will be handled  via Lagrange relaxation.

Since $\hat\gamma_{s,K}(A^\prime,c)>0$, hence, we only consider the case where the elements in $x$ are not all $0$. For any matrix $Q=[q_1,\cdots,q_{n}]\in R^{(m+n)\times n}$ with  $A_2^T Q\leq 0$,
we have
\begin{equation*}
\begin{array}{ll}
Q^TA_1x&=\sum\limits_{\ell=1}^{m+n}\sum\limits_{j=1}^{n}q_{\ell, i} a_{\ell,j}x_j,~i=1,2,\dots,n\\
&=0,
\end{array}
\end{equation*}
where $a_{\ell,j}$ is the element of matrix $A_1$ in the $\ell$-th row and $j$-th column, $q_{\ell,i}$ is the element of matrix $Q$ in the $\ell$-th row and $i-$th column.

Let $C\in R^{n\times n}$ be the diagonal matrix $diag(c_1, c_2, \cdots, c_n)$. For the above $Q$ and by Eq. $\eqref{G}$, we can get

\begin{eqnarray*}
\begin{aligned}
&\hat\gamma_{s,K}(A^\prime,c)\\
&=\max\limits_{\tau, x}\{\sum_{i=1}^{n} \tau^1_ic_ix_i:~\tau\in P_s,~\|x\|_1\leq 1,A_1x=0,~x\geq0\}\\
&\le\max\limits_{\tau\in P_s\atop x\geq0}\left\{\sum_{i=1}^{n} \tau^1_i(c_ix_i-\sum\limits_{\ell=1}^{m+n}\sum\limits_{j=1}^{n}q_{\ell, i} a_{\ell,j}u_j): ~\|x\|_1\leq 1,Q^TA_1 x=0,A_2^T Q\leq 0\right\}\\
&=\max\limits_{\tau\in P_s\atop x\geq0}\left\{\tau^{1T}(Cx-Q^TA_1x): ~ \|x\|_1\leq 1,~Q^TA_1 x=0,A_2^T Q\leq 0\right\}\\
&\le\max\limits_{\tau\in P_s\atop x\geq0}\left\{\tau^{1T}(Cx-Q^TA_1x): ~ \|x\|_1\leq 1,A_2^T Q\leq 0\right\}\\
&=\max\limits_{\tau\in P_s\atop x\geq0}\left\{\tau^{1T}(C-Q^TA_1)x:~\|x\|_1\leq 1,A_2^T Q\leq 0\right\}.
\end{aligned}
\end{eqnarray*}

Hence we can solve the problem
\begin{eqnarray}\label{up2}
\begin{aligned}
\max\limits_{\tau, x}\left\{\tau^{1T}(C-Q^TA_1)x:~x\geq0,\|x\|_1\leq 1,\tau\in P_s,A_2^T Q\leq 0\right\}
\end{aligned}
\end{eqnarray}
to obtain an upper bound of $\hat\gamma_{s,K}(A^\prime,c)$, which is linear in $x$.

In the above problem, the feasible region for $x$ is the convex hull of just $n$ points $ e_i$, $i=1, 2, \cdots, n$, where $e_i$ is the $n$-dimensional vector with the $i$-th component being $1$, and the remaining components being $0$. Therefore, the above problem  can be rewritten as
\begin{equation}\label{upp}
\begin{aligned}
&\max\limits_{\tau, x}\left\{\tau^{1T}(C-Q^TA_1)x:~x\geq0,\|x\|_1\leq 1,\tau\in P_s,A_2^T Q\leq 0\right\}\\
&\leq\max\limits_{\tau, 0< j\leq n}\left\{|\tau^{1T}(C-Q^TA_1)e_j|:~ \tau\in P_s,A_2^T Q\leq 0\right\}\\
&=\max_{0< j\leq n}\left\{\max_{\tau\in P_s}|\tau^{1T}(C-Q^TA_1)e_j|:~A_2^T Q\leq 0\right\}\\
&=\max\limits_{0< j\leq n} \left\{\|(C-Q^TA_1)e_j\|_{s,K,1}:~A_2^T Q\leq 0\right\}.
\end{aligned}
\end{equation}

Define $g_{A_1,C,s,K}(Q)$ as $\max\limits_{0< j\leq n} \|(C-Q^TA_1)e_j\|_{s,K,1}$, and let
\begin{equation*}
\eta_{s,K}(A_1,C,\infty):=\left\{
\begin{array}{cl} \min\limits_{Q} & g_{A_1,C,s,K}(Q)\\
\text{s.t.}&~A_2^T Q\leq 0.
\end{array}
\right.
\end{equation*}
Then
\begin{equation*}
\hat\gamma_{s,K}(A^\prime,c)\leq\eta_{s,K}(A_1,C,\infty).
\end{equation*}
Since $g_{A_1,C,s,K}(Q)$ is easy to compute, so $\eta_{s,K}(A_1,C,\infty)$ is easy to compute.

Further, from \eqref{UP} we have
\begin{equation*}
\begin{aligned}
\max_{\tau\in P_s}\sum_{i=1}^{n} \tau^1_ic_ix_i &\leq\max_{\theta,\nu}\{<A_1x,\theta>+\gamma <x,\nu>: \|\theta\|_*\leq \beta,\|\nu\|_\infty\leq 1\}\\
&=\max_{\theta}\{<A_1x,\theta>:\|\theta\|_*\leq \beta\}+\gamma\|x\|_1\}.
\end{aligned}
\end{equation*}
 Thus, $\gamma$ satisfies that
\begin{equation*}
\begin{aligned}
\max_{\tau\in P_s\atop x\geq0}\left\{\sum_{i=1}^{n} \tau^1_ic_ix_i -\max_{\theta}<A_1x,\theta>: \|\theta\|_*\leq \beta,\|x\|_1\leq1\right\}\leq\gamma.
\end{aligned}
\end{equation*}
Note that $\hat\gamma_{s,K}(A^\prime, c, \beta)$ is the infimum of $\gamma$, hence,
\begin{equation*}
\begin{aligned}
\hat\gamma_{s,K}(A^\prime, c, \beta)=\max\limits_{\tau\in P_s,x\geq0}\left\{\sum_{i=1}^{n} \tau^1_ic_ix_i -\max_{\theta}<A_1x,\theta>:\|\theta\|_*\leq \beta,\|x\|_1\leq1\right\}.
\end{aligned}
\end{equation*}

For any matrix $Q=[q_1,\cdots,q_{n}]\in R^{(m+n)\times n}$ with $A_2^T Q\leq 0$, $\|q_{i}\|_*\leq\beta$ for all $i$, and $Q^TA_1x=0$, we can get
\begin{equation*}
\begin{aligned}
&\hat\gamma_{s,K}(A^\prime,c,\beta)\\
&=\max_{\tau\in P_s\atop {x\geq0 \atop \|x\|_1\leq1}}\left\{\sum_{i=1}^{n} \tau^1_ic_ix_i -\max_{\theta}<A_1x,\theta>:\|\theta\|_*\leq\beta\right\}\\
&\leq\max_{\tau\in P_s\atop {x\geq0 \atop \|x\|_1\leq1}}\left\{\sum_{i=1}^{n} \tau^1_ic_ix_i -<A_1x,q_{i}>: \|q_i\|_*\leq \beta, Q^TA_1x=0,A_2^T Q\leq 0\right\}\\
&=\max_{\tau\in P_s\atop {x\geq0 \atop \|x\|_1\leq1}}\left\{\sum_{i=1}^{n}\tau^1_i(c_ix_i-\sum\limits_{\ell=1}^{m+n}\sum\limits_{j=1}^{n}q_{\ell, i} a_{\ell,j}x_j): \|q_i\|_*\!\leq\! \beta, Q^TA_1x=0,A_2^T Q\leq 0\right\}\\
&=\max_{\tau\in P_s\atop {x\geq0 \atop \|x\|_1\leq1}}\left\{\tau^{1T}(Cx-Q^TA_1x):\! \|q_i\|_*\leq \beta,
Q^TA_1x=0,A_2^T Q\leq 0\right\}\\
&\leq\max_{\tau\in P_s\atop {x\geq0 \atop \|x\|_1\leq1}}\left\{\tau^{1T}(C-Q^TA_1)x: \|q_i\|_*\leq \beta,A_2^T Q\leq 0\right\},
\end{aligned}
\end{equation*}
where $q_i$, $i=1, \dots, n$, is the $i$-th column of matrix $Q$, $a_\ell$, $\ell=1, \dots, n$, is the $\ell$-th column of matrix $A_1$. 

Similar to \eqref{upp}, we can solve the following problem
$$\max_{\tau\in P_s\atop x\geq0} ~\left\{\tau^{1T}(C-Q^TA_1)x: \|x\|_1\leq1,\|q_i\|_*\leq \beta,A_2^T Q\leq 0\right\}.$$
Moreover, it is easy to change the above upper bound of $\hat\gamma_{s,K}(A^\prime,c,\beta)$ to $\eta_{s,K}(A_1,C,\beta)$, which is defined as
\begin{equation}\label{L}
\begin{aligned}
\min\limits_{Q} & \max\limits_{0< j\leq n}~\|(C-Q^TA_1)e_j\|_{s,K,1}\\
\text{s.t.}~ &\|q_i\|_*\leq\beta,~0< i\leq n,\\
&A_2^T Q\leq 0,
\end{aligned}
\end{equation}
where $q_\ell$ is the $\ell$-th column of matrix $Q$.

Obviously,  problem $\eqref{L}$ is a convex programming and is solvable. Similar to the properties of $\hat\gamma_{s,K}(A^\prime,c,\beta)$, $\eta_{s,K}(A_1,C,\beta)$ is a nondecreasing function of $s$, and is a nonincreasing function of $\beta$. Thus we can get an upper bound on $\hat\gamma_{s,K}(A^\prime, c , \beta)$,  by calculating the least upper bound of $g_{A_1,C,s,K}(Q)$ with respect to $Q$, i.e., Eq. \eqref{L}.

%

In addition, according to the definition of $\|\cdot\|_{s,K,1}$, given positive integers $s$ and $t$, we have $\|\cdot\|_{st,K,1}\leq t\|\cdot\|_{s,K,1}$, and
\begin{equation}\label{eta}
\eta_{s,K}(A_1, C, \beta)\leq s\eta_{1,K}(A_1, C, \beta).
\end{equation}
So, the upper bound $\eta_{s,K}(A_1, C, \beta)$ of $\hat{\gamma}_{s,K}(A^\prime, c, \beta)$ can be replaced by $s\eta_{1,K}(A^\prime, C, \beta)$.
This property allows us to reduce the calculation of $\eta_{s,K}(A_1, C)$ to $\eta_{1,K}(A_1, C)$, which greatly reduces the amount of calculation.

 Let $\mathbf{Q}=\{Q: \|q_i\|_*\leq\beta,i=1,\dots,n, ~A_2^T Q\leq 0\}$. According to the definition of $\eta_{s,K}(A_1, C, \beta)$, we have
\begin{equation*}
\begin{array}{lll}
\eta_{1,K}(A_1, C, \beta)&=\min\limits_{Q\in \mathbf{Q}}\max\limits_{0<j\leq n}\|(C-Q^TA_1)e_j\|_{\infty}\\
&=\min\limits_{Q\in \mathbf{Q}}\max\limits_{0<j\leq n}\left\|\left(\begin{array}{ll}|c_1-q_1^T a_1 |,~ \quad |- q_1^ T a_2|, ~&\dots,~ \quad |-q_1^T a_n |\\|-q^T_2a_1|,~\qquad |c_2-q^T_2a_2|, ~&\dots, ~\quad|-q^T_2a_n|\\  &\ddots\\|-q^T_na_1|,~\qquad|-q^T_na_2|, ~&\dots, ~\quad |c_n-q^T_na_n| \end{array}\right)e_j\right\|_{\infty}\\
&=\min\limits_{Q\in \mathbf{Q}}\max\limits_{0<j\leq n}\|(C-A_1^T Q)e_j\|_{\infty}\\
&=\min\limits_{Q\in \mathbf{Q}}\max\limits_{0<j\leq n}\left\|\left(\begin{array}{ll}|c_1- a_1^T q_1 |,~ \quad |- a_1^ T q_2|, ~&\dots,~ \quad |- a_1^T q_n |\\|-a_2^T q_1|,~\qquad |c_2-a_2^T q_2|, ~&\dots, ~\quad|-a_1q^T_n|\\  &\ddots\\|-a_n^Tq_1|,~\qquad|-a_n^Tq_2|, ~&\dots, ~\quad |c_n-a_n^Tq_n| \end{array}\right)e_j\right\|_{\infty}\\
&=\min\limits_{Q\in \mathbf{Q}}\max\limits_{0<j\leq n}\left\|\left(\begin{array}{ll}0\\0\\ \vdots \\c_j\\ \vdots\\0 \end{array}\right)-\left(\begin{array}{ll}a_1^T q_j\\a_2^T q_j \\ \vdots \\a_j^T q_j\\ \vdots \\a_n^T q_j \end{array}\right)\right\|_{\infty}\\
&=\min\limits_{Q\in \mathbf{Q}}\max\limits_{0<j\leq n}\|C_j-A_1^T q_j\|_{\infty},
\end{array}
\end{equation*}
which is equivalent to solve $n$ convex optimization problems of dimension $n$:
\begin{equation}\label{eta1}
\begin{array}{ll}
\eta^j=\min\limits_{q_j\in \mathbf{Q}}\|C_j-A_1^Tq_j\|_{\infty}.
\end{array}
\end{equation}
Obviously, here $\eta_{1,K}(A^\prime, c, \beta)=\max\limits_{0<j\leq n}\eta^j$.

In Theorem \ref{thm8}, the sufficient condition for nonnegative partial $s$-goodness is $\hat\gamma_{s,K}(A^\prime,$ $c)<\frac{1}{2}\min\limits_{0< i\leq n}c_i$. Note that, $\hat\gamma_{s,K}(A^\prime,c)$ takes the value of $\hat\gamma_{s,K}(A^\prime,c,\beta)$ for a large enough $\beta$.
Similarly, $\eta_{s,K}(A_1, C)$ takes the value of $\eta_{s,K}(A_1, C, \beta)$ for the above large enough $\beta$.

Given $A^\prime$, $c$ and $s$, suppose $\hat\gamma_{s,K}(A^\prime,c)<\frac{1}{2}\min\limits_{0< i\leq n}c_i$, such that for every pair of vectors $z^1\in R^n$ and $z^2\in R^{m+n}$, where $z^1\in R^n$ has $s$ nonzero entries, each is equal to 1, there exists a vector $\theta\in R^{m+n}$ such that
\begin{eqnarray*}
\|\big(A_1^{T}\theta\big)-c\circ z^1\|_{\infty}\leq \hat\gamma_{s,K}(A^\prime,c),~
\text{and}~ \big(A_2^{ T}\theta\big)_i\left\{\begin{array}{cc}
=0,&~if~ z^2_i\neq 0;\\
\leq 0,&~if~ z^2_i=0.
\end{array}\right.
\end{eqnarray*}
Next, we need to show that under the assumptions the solution of $\theta$ is bounded, i.e., there exists $\bar{\beta}$ such that  $\|\theta\|_{*}<\bar{\beta}$.

\begin{proposition}\label{pro1}
Suppose $\{u: \|u\|_1\leq \rho, u\in R^{m+n}\}\subseteq
\mathcal{R}(A_1)$, where $\mathcal{R}(A_1)=\{A_1d: \|d\|_1\leq1, d\in R^n\}$. For every $s \leq n$, if $\beta\geq\bar{\beta}=\frac{1}{\rho}(\max\limits_{0< i\leq n}c_i+\frac{1}{2}\min\limits_{0< i\leq n}c_i)$ and $\hat\gamma_{s,K}(A^\prime,c)<\frac{1}{2}\min\limits_{0< i\leq n}c_i$, then $\hat\gamma_{s,K}(A^\prime,c)=\hat\gamma_{s,K}(A^\prime,c,\beta)$.

\end{proposition}
\begin{proof}
 According to the definition of $\hat\gamma_{s,K}(A^\prime,c, \beta)$, we have $\|A_1^T\theta\|_{\infty}<\max\limits_{0< i\leq n}c_i+\frac{1}{2}\min\limits_{0< i\leq n}c_i$. Hence,
\begin{equation*}
\begin{array}{lll}
\max\limits_{0< i\leq n}c_i+\frac{1}{2}\min\limits_{0< i\leq n}c_i&>\|A_1^T\theta\|_{\infty}\\
&=\max\limits_{d}\{d^TA_1^T\theta: \|d\|_1\leq1, d\in R^n\}\\
&=\max\limits_{u}\{u^T\theta: u=A_1d, \|d\|_1\leq1, d\in R^n\}\\
&\geq\max\limits_{u}\{u^T\theta: \|u\|_1\leq\rho\}\\
&=\rho\|\theta\|_{*}.\\
\end{array}
\end{equation*}
Then
$$\|\theta\|_{*}<\frac{1}{\rho}(\max\limits_{0< i\leq n}c_i+\frac{1}{2}\min\limits_{0< i\leq n}c_i),$$
and we define
$$\bar{\beta}=\frac{1}{\rho}(\max\limits_{0< i\leq n}c_i+\frac{1}{2}\min\limits_{0< i\leq n}c_i).$$
The proposition is proven.$\hfill\square$
\end{proof}

According to Proposition \ref{pro1}, if we could  find a lower bound $\bar\beta$ of $\beta$ such that  $\hat\gamma_{s,K}(A^\prime,c, \bar\beta)<\frac{1}{2}\min\limits_{0< i\leq n}c_i$, then  $\hat\gamma_{s,K}(A^\prime,c, \beta)\leq\hat\gamma_{s,K}(A^\prime,c, \bar\beta)<\frac{1}{2}\min\limits_{0< i\leq n}c_i$ for all $\beta\ge\bar{\beta}$, since $\hat\gamma_{s,K}(A^\prime,c,\beta)$ is a nonincreasing function of $\beta$. The same applies to $\eta_{s,K}(A_1, C, \beta)$.

\section{Heuristic algorithm  and examples}\label{chaper5}
In this section, we give three examples to illustrate the proposed non-negative partial $s$-goodness condition for the equivalence between the integer programming problem \eqref{IP} and the weighted linear programming problem \eqref{la3}. To this aim, according to Theorem \ref{thm8}, first we should verify that $(A^\prime,c)$ is nonnegative partially $s$-good. Then, according to  Definition \ref{dingyi1}, the  partially weighted linear programming  problem \eqref{la3} has the unique optimal solution. Further, according to Theorem \ref{thmknownS} or Theorem \ref{thm32}, the optimal solution of problem \eqref{la3} is also the optimal solution of the $l_0$-norm minimization problem \eqref{la2}. Meanwhile, according to Theorem \ref{thm21}, it is also an optimal solution of problem \eqref{IP}.

The main idea of verifying the nonnegative partial $s$-goodness of $(A^\prime,c)$ is as follows. Given the $\bar{\beta}$ in Proposition \ref{pro1}, for $\beta=\bar{\beta}$ and an arbitrary $s$, combining Theorem 8 and the definition of $\eta_{s,K}(A_1, C, \beta)$, it is obvious that $\eta_{s,K}(A_1, C, \bar\beta) < \frac{1}{2}\min\limits_{0<i\leq n}c_i$ is a sufficient  condition for $(A^\prime, c)$ to be nonnegative partially $s$-good.

Let $s^*(A^\prime, c)$ be a lower bound on the largest  $s$ such that $\eta_{s,K}(A_1, C, \bar\beta) < \frac{1}{2}\min\limits_{0<i\leq n}c_i$, which will be obtained  from Eqs. \eqref{eta} and \eqref{eta1}.
In other words, given $s^*(A^\prime, c)$ and $\beta\geq\bar{\beta}$, the value of $\eta_{s^*(A^\prime, c),K}(A_1, C, \beta)$ must be less than $ \frac{1}{2}\min\limits_{0<i\leq n}c_i$.

According to Theorems \ref{thmknownS} and \ref{thm32},  the partially weighted linear programming problem (9) must have a unique optimal solution. However in some cases, the optimal solution of problem \eqref{la3} may not be unique, and we need to adjust $c$, such that problem \eqref{la3} has a unique optimal solution.

To better understand the role of the weight $c$, let us divide a non-negative optimal solution $(x^*,y^*)^T$ with $\|x^*\|_0\leq s$ of problem \eqref{la3} into the following three cases:

$Case~1:$ problem \eqref{la3} has the unique optimal solution $(x^*,y^*)^T$ with $\|x^*\|_0\leq s$;

$Case~2:$ problem \eqref{la3} has multiple optimal solutions $(x^*,y^*)^T$ with $\|x^*\|_0\leq s$, and the vectors $x^*$ have the same sparsity;

$Case~3:$ problem \eqref{la3} has multiple optimal solutions $(x^*,y^*)^T$ with $\|x^*\|_0\leq s$, and the vectors $x^*$ have different sparsity.

Clearly, by Definition \ref{dingyi1}, it is natural that $(A^\prime, c)$ is nonnegative partially $s$-good in $Case~1$. Example \ref{ex1} in this section will illustrate this case.

$Case~2$ shows that there are multiple optimal solutions of problems \eqref{la2} and  \eqref{la3}. At this point, we may adjust $c$ to a suitable value, such that one of the optimal solutions of problem \eqref{la3} is  a unique optimal solution. Then we go to verify that $(A^\prime, c)$ is nonnegative partially $s$-good. Example \ref{ex2} in this section will  illustrate this case.

Note that in $Case~2$, all optimal solutions of problem \eqref{la3} are optimal solutions of problem \eqref{la2}.

$Case~3$ is a very special case, in which problem \eqref{la3} has multiple optimal solutions $(x^*,y^*)^T$, and the vectors $x^*$ have different sparsity. Firstly, suppose that $\eta_{s,K}(A_1, C, \bar\beta)$ $< \frac{1}{2}\min\limits_{0<i\leq n}c_i$, we can obtain $s^*(A^\prime, c)$. Next, we choose one of the optimal solutions of problem \eqref{la3}, which satisfies $\|x^*\|_0\leq s^*(A^\prime, c)$. If such solution cannot be found, then we should find another  $s^*(A^\prime, c)$; otherwise,  we may adjust $c$ to a suitable value, such that this solution is unique, and then we continue to verify that $(A^\prime, c)$ is nonnegative partially $s$-good. Example \ref{ex3} in this section will  illustrate this case.

 The above idea of verifying the nonnegative partial $s$-goodness of $(A^\prime,c)$ can be organized in the following steps. Initially, let $c_i=1,~i=1,2, \dots, n$ and let $\beta=\bar{\beta}=\frac{1}{\rho}(\max\limits_{0< i\leq n}c_i+\frac{1}{2}\min\limits_{0< i\leq n}c_i)$ according to Proposition \ref{pro1}.

$Step~ 1$: According to Eq. \eqref{eta1}, we calculate the value of $\eta_{1,K}(A_1, C, \bar\beta)$. If $\eta_{1,K}(A_1,$ $C, \bar\beta) < \frac{1}{2}\min\limits_{0<i\leq n}c_i$, then go to $Step~ 2$; otherwise (this will happen in $Cases~2$ or $3$), we should use $Step~ 5$ to update $c$, such that $\eta_{1,K}(A_1, C, \bar\beta) < \frac{1}{2}\min\limits_{0<i\leq n}c_i$, and go to $Step~ 2$.

$Step~ 2$: Since $s$ is not known at present, we suppose $\eta_{s,K}(A_1, C, \bar\beta) < \frac{1}{2}\min\limits_{0<i\leq n}c_i$. Then according to \eqref{eta},
$$s^*(A^\prime, c)=\lfloor\frac{\frac{1}{2}\min\limits_{0<i\leq n}c_i}{\eta_{1,K}(A_1, C, \bar\beta)}\rfloor.$$

$Step~ 3$: Consider an optimal solution $(x^*,y^*)$ of problem \eqref{la3}. If the solution of problem \eqref{la3} is  unique and $\|x^*\|_0=s$, then we compare $s$ with  $s^*(A^\prime, c)$. If $s=s^*(A^\prime, c)$, then we verify whether $s^*(A^\prime, c)\eta_{1,K}(A_1, C,$ $ \bar\beta) < \frac{1}{2}\min\limits_{0<i\leq n}c_i$. When it holds, go to $Step~ 4$.

Otherwise, such as $Case~3$, not all solutions satisfy $s\eta_{1,K}(A_1, C, \bar\beta) < \frac{1}{2}\min\limits_{0<i\leq n}c_i$. So we choose a solution with $\|x^*\|_0=s^*(A^\prime, c)$, and use $Step~ 5$ to update $c$, such that this optimal solution of problem \eqref{la3} is the unique optimal solution. Next,  we verify that whether $s^*(A^\prime, c)\eta_{1,K}(A_1, C, \bar\beta) < \frac{1}{2}\min\limits_{0<i\leq n}c_i$. When it holds, go to $Step~ 4$; otherwise, update $c$ again.

$Step~ 4$: According Theorem \ref{thm8}, it implies that $(A^\prime, c)$ is nonnegative partially $s$-good. Stop the algorithm.




$Step~ 5$: Update $c$ as follows: for the maximum component $x_i^*$ in $Step~ 3$, select $c_i$ such that $0<c_i\leq \bar\beta$. For the minimum component $x_j^*$ in $Step~ 3$, select $c_j$ from $(\bar\beta, \frac{3}{2}\bar\beta)$. The other components in $c$ are randomly selected from $[c_i, c_j]$, i.e., $c_i$ and $c_j$ are the minimum and maximum components of $c$ respectively. Let $\bar{\beta}=\frac{1}{\rho}(\max\limits_{0< i\leq n}c_i+\frac{1}{2}\min\limits_{0< i\leq n}c_i)$.

Below are three examples we provide.  Example \ref{ex1} does not comply with $TUM$, Example \ref{ex2} does not comply with $TDI$, and Example \ref{ex3} is neither $TUM$ nor $TDI$ compliant.

\begin{example}\label{ex1}
Let
\begin{eqnarray*}
A=\left(\begin{array}{lll}
1 & 2 &~ 0 \\
0 & 1 &~ 1 \\
1 & 0 & 2
\end{array}\right),  \quad b=\left(\begin{array}{l}
1 \\
1 \\
1
\end{array}\right).
\end{eqnarray*}

For $c=1$, given $\bar\beta=0.563$, according to Eq. \eqref{eta1}, we can calculate the value of $\eta_{1,K}(A_1,C, \bar\beta)$ is $0.2188$. Suppose $\eta_{s,K}(A_1, C, \bar\beta) < \frac{1}{2}\min\limits_{0<i\leq n}c_i$, then according to \eqref{eta}, $s^*(A^\prime, c)=\lfloor\frac{0.5}{0.2188}\rfloor=2$. For $c=1$, problem \eqref{la3} has an optimal solution  $\left(x^{1}, y^{1}\right)^T=\left(\left(0, \frac{1}{2}, \frac{1}{2}\right),\left(0,0,0, 1,\frac{1}{2},\frac{1}{2}\right)\right)^{T}$ and $\|x^1\|_0\leq2$. According to \eqref{eta}, $\eta_{s,K}(A_1, C, \bar\beta)\leq2\eta_{1,K}$ $(A_1,C, \bar\beta)$ $=0.4376 < \frac{1}{2}$. Then, by Theorem \ref{thm8}, $(A^\prime,c)$ is nonnegative partially $s$-good. Hence, according to Definition \ref{dingyi1}, $(x^1,y^1)^T$ is the unique optimal solution of  problem \eqref{la3}.  Next, according to Theorem \ref{thm32}, $(x^1,y^1)^T$ is an optimal solution of problem \eqref{la2}. So by Theorem \ref{thm21}, $(0,1,1)^T$ is an optimal solution of  problems \eqref{IP} and \eqref{la1} respectively.
\end{example}

It must be remarked that, sometimes problem (\ref{la3}) has multiple nonnegative optimal solutions  when $c=1$, and the vectors $x$ have the same sparsity. This situation is contradict to Definition \ref{dingyi1}. So the coefficient $c$ must be adjusted, such that  problem (\ref{la3}) has a unique optimal solution. The following Example \ref{ex2} shows that this can be achieved.

 \begin{example}\label{ex2}
Let
\begin{eqnarray*}
A=\left(\begin{array}{lll}
1 & 0 & 0 \\
1 & 1 & 0 \\
0 & 1 & 1
\end{array}\right),~ b=\left(\begin{array}{l}
0 \\
\frac{3}{2} \\
\frac{1}{2}
\end{array}\right).
\end{eqnarray*}

For $c=1$, given $\bar\beta=0.5$, according to Eq. \eqref{eta1}, we can calculate the value of $\eta_{1,K}(A_1,C, \bar\beta)$ is $0.5$. For $c=1$, problem \eqref{la3} has two optimal solutions, which are $\left(x^{1}, y^{1}\right)^T=\left(\left(1, \frac{1}{2}, 0\right),\left(1,0,0,0,\frac{1}{2},1\right)\right)^{T}$ and $\left(x^{2}, y^{2}\right)^T=\left(\left(\frac{3}{4}, \frac{3}{4}, 0\right), \left(\frac{3}{4}, 0, \frac{1}{4}, \frac{1}{4}, \frac{1}{4},1\right)\right)^{T}$,  $c^Tx^1=c^Tx^2$, and both vectors $x^1$ and $x^2$ are $2$-sparse.


Next, let $c=(0.5,0.7,0.8)$. Given $\bar\beta=0.7$, we have $\eta_{1,K}(A_1,C, \bar\beta)=0.1$. Suppose $\eta_{s,K}(A_1, C, \bar\beta) < \frac{1}{2}\min\limits_{0<i\leq n}c_i$, then according to \eqref{eta}, $s^*(A^\prime, c)=\lfloor\frac{0.25}{0.125}\rfloor=2$. For $c=(0.5,0.7,0.8)$, 
 problem \eqref{la3} has an optimal solution
$\left(x^{1}, y^{1}\right)^T=\left(\left(1, \frac{1}{2}, 0\right), \left(1, 0, 0, \right.\right.$ $ \left.\left.0,\frac{1}{2},1\right)\right)^{T}$, and $\|x^*\|_0\leq2$. According to \eqref{eta}, $\eta_{s,K}(A_1, C, \bar\beta)\leq2\eta_{1,K}(A_1,C)=0.2 < \frac{1}{2}\min\limits_{0<i\leq n}c_i=0.25$. Then, by Theorem \ref{thm8}, $(A^\prime,c)$ is nonnegative partially $s$-good. Hence, according to Definition \ref{dingyi1}, $(x^1,y^1)^T$ is the unique optimal solution of problem \eqref{la3}. Furthermore, according to Theorem \ref{thmknownS}, $(x^1,y^1)^T$ is an optimal solution of problem \eqref{la2}. Since $\|x^1\|_0=\|x^2\|_0$, it is natural that $(x^2,y^2)^T$ is also an optimal solution to problem \eqref{la2}. So by Theorem \ref{thm21}, $(1, 1, 0)$ is an optimal solution of  problems \eqref{IP} and \eqref{la1}  respectively.
\end{example}

Different from Example \ref{ex2}, another situation to be aware of is that, problem (\ref{la3}) has multiple nonnegative optimal solutions, and the vectors $x$ have different sparsities, when $c=1$.  This situation is also contradict to Definition \ref{dingyi1}. So the coefficient $c$ must be adjusted too, such that problem (\ref{la3}) has a unique optimal solution. The following Example \ref{ex3} shows that this can be achieved. 

\begin{example}\label{ex3}
Let
\begin{eqnarray*}
A=\left(\begin{array}{lll}
1 & 2 &~ 0 \\
0 & 1 &~ 1 \\
2 & 0 & 1
\end{array}\right),  \quad b=\left(\begin{array}{l}
0 \\
\frac{1}{2} \\
\frac{1}{3}
\end{array}\right).
\end{eqnarray*}

For $c=1$, given $\bar\beta=0.375$, according to Eq. \eqref{eta1}, we can calculate the value of $\eta_{1,K}(A_1,C, \bar\beta)$ is $0.2917$. Suppose $\eta_{s,K}(A_1, C, \bar\beta) < \frac{1}{2}\min\limits_{0<i\leq n}c_i$, then according to \eqref{eta}, $s^*(A^\prime, c)=\lfloor\frac{0.5}{0.2917}\rfloor=1$. For $c=1$, problem \eqref{la3} has two optimal solutions, they are: $\left(x^{1}, y^{1}\right)^T=\left(\left(0, \frac{1}{6}, \frac{1}{3}\right), \left(\frac{1}{3},0,0,1,\frac{5}{6},\frac{2}{3}\right)\right)^{T}$ with $\|x^1\|_0=2$, and $\left(x^{2}, y^{2}\right)^T=\left(\left(0, 0, \frac{1}{2}\right),(0, 0,\frac{1}{6}, 1, 1, \frac{1}{2})\right)^{T}$ with $\|x^2\|_0=1$, $c^Tx^1=c^Tx^2$. Clearly, $\|x^1\|_0>s^*(A^\prime, c)$, $\|x^2\|_0=s^*(A^\prime, c)$.

Next, let $c=(0.5,0.35,0.3)$. Given $\bar\beta=0.7$, we have $\eta_{1,K}(A_1,C, \bar\beta)=0.1$. Suppose $\eta_{s,K}(A_1, C, \bar\beta) < \frac{1}{2}\min\limits_{0<i\leq n}c_i$, then according to \eqref{eta}, $s^*(A^\prime, c)=\lfloor\frac{0.15}{0.1}\rfloor=1$. For $c=(0.5,0.7,0.8)$, problem \eqref{la3} has an optimal solution $\left(x^{2}, y^{2}\right)^T=\left(\left(0, 0, \frac{1}{2}\right),\left(0, 0,\frac{1}{6}, \right.\right. $\\$\left.\left.1, 1, \frac{1}{2}\right)\right)^{T}$ with $\|x^2\|_0$ $=1$. According to \eqref{eta}, $\eta_{s,K}(A_1, C, \bar\beta)\leq\eta_{1,K}(A_1,C, \bar\beta)=0.1 < \frac{1}{2}\min\limits_{0<i\leq n}c_i=0.15$. Hence, by Theorem \ref{thm8}, $(A^\prime,c)$ is nonnegative partially $s$-good, and according to Definition \ref{dingyi1}, $(x^2,y^2)^T$ is the unique optimal solution of  problem \eqref{la3}.  Furthermore, according to Theorem \ref{thm32}, $(x^2,y^2)^T$ is an optimal solution of problem \eqref{la2}. So by Theorem \ref{thm21}, $(0,0,1)^T$ is an optimal solution of  problems \eqref{IP} and \eqref{la1} respectively.
\end{example}

\section{Conclusion}\label{chaper6}

In this paper, we studied the equivalence of a 0-1 linear program to a weighted linear programming problem. Firstly, we prove the equivalence between the integer programming problem and a sparse minimization problem. Next, we define the nonnegative partial $s$-goodness of the constraint matrix and the weight vector in the objective function of the weighted linear programming problem. Utilizing two quantities $\gamma_{s,K}(\cdot)$ and $\hat{\gamma}_{s,K}(\cdot)$ of the nonnegative partial $s$-goodness, we propose a necessary  and a sufficient condition for the constraint matrix and weighted vector to be nonnegative partial $s$-good. Since it is difficult to calculate the two quantities, we further provide an efficiently computable upper bound of $\hat{\gamma}_{s, K}(A^\prime, c, \beta)$, such that the above sufficient condition is verifiable. It is worthy of mentioning that the objective coefficient $c$ of the weighted linear programming problem is not fixed. When the weighted linear programming problem has multiple optimal solutions, we may adjust $c$ so that the weighted linear programming problem has only a unique optimal solution. At the end, we provide three examples to illustrate the theory in this article.


%
%

\bibliographystyle{spmpsci}      
\bibliography{mybib}   

\end{document}